\documentclass[10pt]{amsart}

\usepackage{amscd, amsmath, amssymb, latexsym}
\usepackage[dvips]{graphics}
\usepackage{epsfig, psfrag}

\newcommand{\refT}[1]{Theorem~\ref{T:#1}}
\newcommand{\refC}[1]{Corollary~\ref{C:#1}}
\newcommand{\refP}[1]{Proposition~\ref{P:#1}}
\newcommand{\refD}[1]{Definition~\ref{D:#1}}
\newcommand{\refL}[1]{Lemma~\ref{L:#1}}

\newcommand{\Z}{{\mathbb Z}}
\newcommand{\R}{{\mathbb R}}

\newcommand{\Sto}{{\Psi^o_{n+2}}}
\newcommand{\Ston}{{\Psi^o_{n+1}}}

\newcommand{\C}{{\mathcal C}}
\newcommand{\D}{{\mathcal D}}

\newcommand{\F}{{\mathcal F}}
\newcommand{\G}{{\mathcal G}}
\newcommand{\Ha}{{\mathcal H}}
\newcommand{\X}{{\mathcal X}}

\newcommand{\E}{{{\rm Emb}(I, M)}}
\newcommand{\Int}{{\rm Int}}

\newcommand{\Em}{{\mathcal E}}
\newcommand{\wt}{\widetilde}

\newcommand{\und}{\underline}
\newcommand{\la}{\langle [}
\newcommand{\ra}{] \rangle }
\newcommand{\pn}{\mathcal{P}(\underline{n})}
\newcommand{\pno}{\mathcal{P}(\underline{n+1})}
\newcommand{\pnn}{\mathcal{P}_\nu(\n)}
\newcommand{\pnno}{\mathcal{P}_\nu(\underline{n+1})}
\newcommand{\n}{\underline{n}}
\newcommand{\no}{\underline{n+1}}

\newcommand{\skiphome}[1]{#1}
\newcommand{\iop}{{{\mathcal S}{\it iop\/}}}

\newtheorem{theorem}{Theorem}[section]
\newtheorem{lemma}[theorem]{Lemma}
\newtheorem{proposition}[theorem]{Proposition}
\newtheorem{corollary}[theorem]{Corollary}

\newtheorem{definition}[theorem]{Definition}

\theoremstyle{remark}

\begin{document}

\title{The topology of spaces of knots: cosimplicial models}
\author{Dev P. Sinha}
\address{Department of Mathematics, University of Oregon, Eugene OR}


\begin{abstract}
We present two models for the space of knots which have endpoints at fixed boundary points in a manifold with boundary, one model defined as an inverse limit of  spaces of maps between configuration spaces and another which is cosimplicial.   These models build on the calculus of isotopy functors and  are weakly homotopy equivalent to knot spaces when the ambient dimension is greater than three. The mapping space model, and the evaluation map on which it builds,  is suitable for analysis through differential topology. The cosimplicial model gives rise to  spectral sequences which converge to cohomology and homotopy groups of spaces of knots when they are connected.    We explicitly identify and establish vanishing lines in these spectral sequences.

 \end{abstract}


\maketitle

{\tableofcontents}

\section{Introduction}

We study the topology of $\E$,  the space of smooth
embeddings of the interval into $M$ where the embeddings are constrained to map the
boundary of the interval to fixed boundary points of $M$, with fixed tangent
vectors at those boundary points.   When $M = I^N$ this space is homotopy equivalent to the space
of  long knots.

We build on the approach of Goodwillie and his collaborators \cite{GK,
GKW, GW, We}, which is known as the calculus of embeddings.  From their approach
we produce two equivalent models.  The first is the mapping space model, which provides
the proper setting for an evaluation map suitable for geometric study \cite{BCSS02}.  
The second model  is
cosimplicial, analogous to the cosimiplicial model for loop spaces and especially
convenient for calculations of homotopy and cohomology groups.  Both of
these models use completions of configuration spaces,
due to Fulton and MacPherson \cite{FM} and done in the category
of manifolds by Axelrod and Singer \cite{AxelSing}, but with some 
changes needed for the cosimplicial model \cite{Sinh04}. 

The cosimplicial model gives rise to spectral sequences for both cohomology and
homotopy groups of $\E$ which converge when $M$ is simply connected
and has dimension four or greater, so that this knot space is connected.
We explicitly identify these spectral sequences and establish vanishing lines.  
When $M$ is $I^N$, the
cohomology spectral sequence is reminiscent of those of Vassiliev
\cite{Vass, Vass97} and Kontsevich \cite{Kont94, Kont00}, as further studied by
Tourtchine \cite{To00, Tour04}.   Indeed, Tourtchine has announced comparison theorems
showing that these all have isomorphic $E_2$ pages.
The homotopy spectral sequence  gives the first
computations of homotopy groups of embedding spaces.  The 
$E^1$-term was conjectured rationally by Kontsevich \cite{Kont00}. 
We study the rational spectral sequence for knots in even-dimensional Euclidean
spaces in \cite{2}, giving explicit computations
in low dimensions.   Tourtchine and Lambrechts have recast this
spectral sequence in terms of a new form of graph homology \cite{LaTo06}.

When the dimension of $M$ is three it is not known whether our models
are equivalent to corresponding knot spaces.   We can still pull back
zero-dimensional cohomology classes (that is, knot invariants) from our
models.  These invariants are connected in various ways with the
theory of finite-type invariants \cite{BCSS02}, and in particular
are at least as strong as rational finite-type invariant theory 
by the  work of Volic \cite{Voli03}.  
In further work, we hope to extend the constructions of \cite{BCSS02} to
higher degrees in order to define finite-type knot invariants through
differential topology.

The study of spaces of knots is highly active.   
In \cite{Sinh03} we build on the models of this paper to connect spaces of knots
with the theory of operads, resolving conjectures
of Kontsevich from \cite{Kont00} which were inspired by Tourtchine's
remarkable algebra \cite{To00}.  Vassiliev uses analogies between knot
spaces and hyperplane arrangements to make further progress on realizing
classes in his spectral sequence \cite{Vass03}.  Arone, Lambrechts, Tourtchine and Volic
have established collapse of the rational
homotopy spectral sequence of this paper for knots in Euclidean space using 
coformality of the Fulton-MacPherson operad \cite{ALTV07}.  They have also 
announced collapse of the cohomology spectral sequence, as conjectured by Kontsevich and 
Vassiliev.   Budney has shown that spaces of framed long knots in
Euclidean space have a little two-cubes action, as first conjectured by Tourtchine.
For knots in three dimensions, Budney has used techniques of Hatcher \cite{Hatc01} and 
decomposition
theorems for three manifolds to show that this two-cubes action is free \cite{Budn03}
and identify the homotopy type of the knot space \cite{Budn05}.
We hope to survey these developments at another time.

\subsection{Basic definitions, notation and conventions}

We choose a variant of $\E$ so that a knot depends only on its image.
To do so requires a Riemannian metric on $M$.  Some of our constructions will
depend on this metric, but changing the metric will
always result in changing spaces in question by a homeomorphism.  Thus
we make little mention of the metric in general.

\begin{definition}
Let $\E$  be the space of injective $C^1$ maps of constant
speed from $I$ to $M$,  whose values and unit tangent vectors at $0$ and
$1$ are specified by fixed inward- (respectively outward-) pointing
unit tangent vectors in $\partial M$, with the $C^1$-Whitney topology.
\end{definition}

We use spaces of ordered configurations, and modifications thereof,
extensively.  There are at least two lines of notation for these spaces in the
literature.  We prefer the notation $C_n(M)$ rather than $F(M,n)$.  
Note, however, that $C_n(M)$ in this paper is
$C^0_n(M)$ in \cite{BT} and that $C_n[M]$ in this paper is $C_n(M)$ in
\cite{BT}.   Indeed, we must pay close attention to the
parentheses in our notation, because for the sake of brevity of notation 
they account for all of the distinction between the variants configuration
spaces:  $C_n(M)$ is the standard ``open'' configuration space;
$C_n[M]$ is the Fulton-MacPherson completion, the canonical
closure; $C_n\la M \ra$ is a quotient of $C_n[M]$.

 Let $\n$ denote the set $\{1, \ldots, n\}$, our most
common indexing set, and let $C_m(S)$ denote the set of distinct
$m$-tupes in an ordered set $S$.  Let $[n] = \{0, 1, \ldots, n\}$.

To facilitate working with products, set  $X^S= {\rm Maps}(S, X)$ for
any finite set $S$.  Consistent with this, if $\{X_s\}$ is a
collection of spaces indexed by $S$, we let $(X_s)^S = \prod_{s \in S}
X_s$.  For coordinates in either case we use
$(x_s)_{s\in S}$ or just $(x_s)$ when $S$ is understood. Similarly,  a
product of maps $\prod_{s\in S} f_s$ may be written $(f_s)_{s \in S}$ or
just $(f_s)$.

We work with stratifications of spaces in a fairly naive sense.  A
stratification of a space $X$ is a collection of disjoint subspaces $\{X_c\}$ such
that the intersection of the closures of any two strata is the closure
of some stratum.   We associate a poset to a
stratification by saying that $X_c \leq X_d$ if $X_c$ is in the closure
of $X_d$.  If two spaces
$X$ and $Y$ are stratified with the same associated posets, 
then a stratum-preserving map is one in which the closure of $X_c$ maps
to the closure of $Y_c$ for every $c$.

We use limits and homotopy limits of diagrams extensively.  We recommend
\cite{GoII} for a more complete introduction.
A diagram of spaces is a functor from a
small category $C$, which we sometimes refer to as an indexing category, 
to the category of spaces.  Denote the
realization of the nerve of $C$ by $|C|$.  
If $c$ is an object of $C$, the category $C
\downarrow c$  has objects which are maps with target 
$c$ and morphisms given by morphisms in $C$ which commute with these 
structure maps.  Note that $|C \downarrow c|$ is contractible 
since $C \downarrow c$ has a final object, namely $c$  mapping to itself by the 
identity morphism.  A morphism from $c$ to $d$ induces a map  
from $C \downarrow c$ to $C \downarrow d$, so that $|C\downarrow -|$ 
is a functor from $C$ to spaces. 

\begin{definition}
The homotopy limit of a functor $\F$ from a small category $C$ to the 
category of spaces is ${\rm Nat}(|C \downarrow -|, \F)$,
the space of natural transformations from $|C \downarrow -|$ to
$\F$.
\end{definition}

Let $C$ be a  poset with a unique maximal element $m$, so that 
$C \cong C \downarrow m$, and with unique greatest lower bounds.  
For any object 
$c$ of $C$ the category $C \downarrow c$ is naturally a sub-category of $C 
\downarrow m$. 
The $|C \downarrow c|$ define a stratification of $|C| = |C \downarrow 
m|$, where the ``open'' strata are $|C \downarrow c| - \bigcup_{d < c} |C \downarrow c|$.
The poset associated to this stratification is $C$.  

Let $\F$ be a functor from 
$C$ to spaces in which all morphisms are inclusions of NDR pairs 
in which the subspace in each pair is closed.
Then $\F(m)$ is stratified by the $\F(c)$ with associated poset 
isomorphic to $C$ once again.  

\begin{proposition}\label{P:inclimit}
Let $C$ and  $\F$ be as above.  
The homotopy limit of $\F$ is the space of all stratum-preserving 
maps from $|C \downarrow m| = |C|$ to $\F(m)$. 
\end{proposition} 

Some basic properties of homotopy limit are the following  
(see \cite{GoII} and \cite{BK}).

\begin{proposition}\label{P:holim}
\begin{enumerate} 
\item The homotopy limit is functorial. 
\item If there 
is a natural transformation $\F \to \G$ in which each 
$\F(c) \to \G(c)$ is a (weak) homotopy equivalence, then the induced map 
on homotopy limits is a (weak) homotopy equivalence. 
\item There is a canonical map from the limit to the homotopy limit.  
\item (``Quillen's Theorem A'') If a functor $\F$ from $D$ to $C$ is left 
cofinal then composition with $\F$ induces a weak equivalence on 
homotopy limits.  For $c \in C$, let $\F \downarrow c$ be the category whose objects are 
pairs $(d, f)$ where 
$d$ is an object of $D$ and $f$ is a morphism from $\F(d)$ to $c$.
Morphisms are given by morphisms $g$ in $D$ such that $\F(g)$ commutes
with the structure maps.  A functor $\F$ is left cofinal if $|\F 
\downarrow c|$ is contractible for every object $c \in C$.  
\end{enumerate} 
\end{proposition}

\subsection{Acknowledgements}

The author is deeply indebted to Tom Goodwillie.  As mentioned in
\cite{GKW} (and also briefly in \cite{BT}), Goodwillie has known for
some time that one should be able to use Theorem~\ref{T:cusimp} to
give a cosimplicial model and thus resulting spectral sequences.
Thanks
also go to Victor Tourtchine for many helpful comments.

\section{Heuristic understanding of the mapping space and cosimplicial
models}\label{S:heur}

Given a knot $\theta \colon I \to M$ and $n$ distinct points on the unit interval, one may
produce $n$ distinct points in the target manifold by evaluating the knot
at those points.  Because a knot has a nowhere vanishing derivative 
(``no infinitesimal self-intersections''), one may
in fact produce a collection of $n$ unit tangent vectors at these distinct
points.   Let
$\Int(\Delta^n)$ be the open $n$-simplex and let $C'_n(M)$ be defined by
the pull-back square
$$
\begin{CD}
C'_n(M)@>>> (STM)^n \\ 
@VVV   @VVV \\
C_n(M) @>>> M^n,
\end{CD}
$$
where $STM$ is the unit tangent bundle of $M$.

We now define the evaluation map as follows.

\begin{definition}
Given a knot $\theta$ let
$ev_n(\theta) \colon \Int(\Delta^n) \to C'_n(M)$ be defined by
$$(ev_n(\theta))(t_1, \cdots t_n) = (u(\theta'(t_1)), \cdots,
u(\theta'(t_n))),$$ where $u(v)$ is the unit tangent vector in the
direction of the tangent vector $v$. Let $ev_n \colon \E
\to {\text{Maps}}(\Int(\Delta^n), C'_n(M))$ be the map which sends
$\theta$ to $ev_\theta$.
\end{definition}

The evaluation map is sometimes called a Gauss map, as classically
it is used to define the linking number.  We cannot  expect to use it
to study the homotopy type of the knot spaces as it stands since
$ev_n$ agrees in the homotopy category with the map $\E \to C'(M)$ which evaluates
a knot at a single point.  In order to define a mapping space which could possibly reflect the
topology of the embedding space, we need to add boundaries to configuration spaces
and impose boundary conditions. 
These boundaries are part of Fulton-MacPherson completions,
which have the same homotopy type as the open 
configuration spaces and  are functorial for embeddings.  The 
Fulton-MacPherson completion is a manifold with corners, and for 
any $\theta$ the evaluation map 
$ev_n(\theta)$ respects the resulting stratification.  Our mapping space models are 
essentially the spaces of stratum-preserving maps from the space of configurations
in the interval to the space of configurations in $M$ (with one important
additional technical condition).  The maps from 
$\E$ to these models are extensions of the $ev_n$.

To motivate our cosimplicial model, recall 
the cosimplicial model for the based loop space, $\Omega M$ (see 
for example \cite{Re70}).  The $n$th entry of this cosimplicial model is 
given by the Cartesian product $M^{\n}$.  The coface maps are  diagonal 
(or ``doubling'') maps, and the codegeneracy maps are projections (or
``forgetting'' maps). The map from the loop space to the $n$th total space 
of this cosimplicial space is the adjoint of an evaluation map from the 
simplex to $(M)^{\n}$, and it is a homeomorphism 
if $n \geq 2$.  

To make cosimplicial models of knot spaces we try 
to replace the cosimplicial entry $(M)^{\n}$ by the configuration space $C'_n(M)$, 
since  $C'_n(M)$ is a natural 
target for the evaluation map for embeddings.    
The codegeneracy maps can be defined as for the loop space, by 
forgetting a point in a configuration.  The coface maps 
are problematic, as we cannot double a point in a configuration 
to get a new configuration of distinct points.  We are tempted to add a point close to 
the point which needs to be doubled, but in order for the composition
of doubling and forgetting to be the identity, we need to add a 
point which is ``infinitesimally close.''  The appropriate technical 
idea needed to overcome this difficulty is once again that of the 
Fulton-MacPherson completion.  But while 
this completion has diagonal maps, these maps do not satisfy the cosimplicial axioms.  
We thus use a variant of this completion 
which admits a cosimplicial structure.

\section{Goodwillie's cutting method for knot spaces}\label{S:cutting}

One version of  Goodwillie's cutting method  approximates the space of 
embeddings of a manifold $M$ using
embeddings of $M - A$, for codimension zero submanifolds $A$ of $M$. 
We fix a collection  of disjoint closed sub-intervals of $I$ by setting 
$J_i = (\frac{1}{2^i}, \frac{1}{2^i} + \frac{1}{10^i})$.

\begin{definition} 
Define $E_S(M)$ for $S \subseteq 
\und{n} = \{1, \cdots, n\}$ to be the space of embeddings of $I - \bigcup_{s \in S} 
J_s$ in $M$ whose speed is constant on each component, topologized with
with $C^1$-topology.   
\end{definition} 

By convention, an embedding of any interval containing $0$ or $1$
must send those points to the designated points on the boundary
of $M$, with the designated unit tangent vectors at those points.
Because the speed is constant on each component, an element of
$E_S(M)$ is determined by its image in $M$.

If $S \subseteq S'$ there is a restriction map from $E_S(M)$ to 
$E_{S'}(M)$.  These restriction maps  commute, so to a 
knot we can associate a family of compatible elements of 
$E_S$ for every non-empty $S$.  Conversely, such a compatible family  
 determines a knot if $n >2$.   Motivated by notions of ``higher-order excision'',
Goodwillie's cutting method uses 
families of punctured knots compatible only up to isotopy to approximate
the space of knots.  Such families of punctured knots are described through
homotopy limits.   
Recall from \cite{GoII} the language of cubical diagrams.

\begin{definition} 
Let $\pn$ be the category of all subsets of $ \underline{n}$ 
where morphisms are defined by inclusion.  Let $\pnn$ be
the full subcategory of non-empty subsets.  A cubical (respectively sub-cubical)
diagram of spaces is a functor 
from $\pn$ (respectively $\pnn$) 
to the category of based spaces. 
\end{definition} 

The nerve of $\pn$ is an $n$-dimensional  
cube divided into simplices.   The nerve of 
$\pnn$ consists of $n$ faces of that cube and is isomorphic 
to the barycentric subdivision of the $(n-1)$-simplex. 

\begin{definition}\label{D:ej} 
 Let $\overline{\Em_n}(M)$ be the cubical diagram which sends $S \in 
\pno$ to the space of embeddings $E_S(M)$ and sends the inclusion of $S \subset S'$ to 
the appropriate restriction map.  Let $\Em_n(M)$ be the restriction of $\overline{\Em_n}(M)$ 
to $\pnno$.
\end{definition}  

A cubical diagram such as $\overline{\Em_n}(M)$ determines a map from the 
initial space in the cube, in this case $\E$, to the homotopy limit of 
the rest of the cube. 

\begin{definition}\label{D:defpn}
Let $P_n \E$ be the homotopy limit of $\Em_{n}(M)$.  Let $\alpha_n$ be 
the canonical  map from $\E$, which is the initial space in 
$\overline{\Em_{n}}(M)$, to $P_n \E$.  
\end{definition} 

The space $P_n \E$ is a degree $n$ polynomial approximation to the space of 
knots in the sense of the calculus of isotopy functors \cite{We}.    In particular,
because removing a single interval $J_{i}$ results in a contractible embedding
space,  $P_{1} \E \simeq \Omega {\rm Emb}\left(I - (J_{1} \cup J_{2}), M\right)$.  From 
Proposition~\ref{P:incl} below it follows that
$ {\rm Emb}\left(I - (J_{1} \cup J_{2}), M\right)$ is 
homotopy equivalent to $\Omega STM$, the loop space of the unit tangent
bundle of $M$, which is known in turn to be homotopy equivalent to the space
of immersions ${\rm Imm}(I, M)$.
This space is the ``linear'' approximation to the space of embeddings,
because immersions exhibit a Mayer-Vietoris property.

\begin{theorem}[\cite{GK}] \label{T:good} 
Let $M$ be a manifold whose dimension 
is four or greater.  The map $\alpha_n$ 
from $\E$ to $P_n \E$  is $(n-1)({\rm dim} \; M -3)$-connected. 
\end{theorem} 

In Section~1.A of
\cite{Good06},  Goodwillie provides 
remarkably elementary arguments using dimension counting and the 
Blakers-Massey theorem which prove a version of  this theorem 
with weaker connectivity estimates for 
$M$ of dimension five and higher.   To get the best connectivity estimates
and to apply in dimension four
one needs the deep arguments of Goodwillie and Klein \cite{GK}.

The category $\pnn$ is a sub-category of
$\pnno$, through the standard inclusion of $\{1, \cdots, n\}$ into 
$\{1,\cdots, n+1\}$ for definiteness.  
Because our choices of the $J_{i}$ in the definition of $\Em_{n}(M)$ 
 are compatible, $P_{n} \E$ maps to 
$P_{n-1} \E$ through a restriction map $r_n$ such that $r_n \circ \alpha_{n} 
= \alpha_{n-1}$.  By \refT{good}, the maps $\alpha_n$ induce
isomorphisms on homology and homotopy groups through a range which 
increases linearly in $n$, so we deduce the following.

\begin{corollary}[\cite{GK}]\label{C:weq} 
If the dimension of the ambient manifold is greater than three, the 
$\alpha_{n}$ give rise to a map 
from $\E$ to the homotopy inverse limit of   
$$P_0 \E \leftarrow P_1 \E \leftarrow P_2 \E \leftarrow \cdots$$ 
which  is a weak equivalence. 
\end{corollary}  

Because the spaces of punctured knots
$E_S(M)$ are  homotopy equivalent to configuration spaces 
$C'_n(M)$ (see the proof of
Proposition~\ref{P:incl}), we are led to search for geometric models equivalent to 
$P_n \E$ which involve these spaces.  
Completions of configuration spaces are essential 
to our construction of such models. 

\section{Fulton-MacPherson completions}\label{S:FMtrees} 

We use two different versions of completions of
configuration spaces.  For varieties these were defined in the seminal paper \cite{FM}.
In the setting of manifolds they were first 
defined in \cite{AxelSing} and further developed in \cite{Mark99, Gaif03, Sinh04}.  In this
section we review definitions and results needed from \cite{Sinh04},
where full statements and proofs  are provided.

\subsection{Basic definitions and properties}

To define these completions, we need to fix an isometric embedding of our manifold $M$
in some Euclidean space $\R^{N+1}$, which is the identity if $M$ is $\R^{N+1}$ 
or the standard inclusion if $M = I^{N+1}$.
All constructions which we make are ultimately, up to homeomorphism, 
independent of this ambient embedding.

\begin{definition}\label{D:basiccomp}
\begin{enumerate}

\item If $M$ is a smooth manifold let $C_n(M)$ be the subspace of  
$(x_i) \in M^{\n}$ such that $x_i \neq x_j$ if $i \neq j$.  Let
$\iota$ denote the inclusion of $C_n(M)$ in $M^{\n}$.   

\item For $(i,j) \in C_2(\n)$,
let $\pi_{ij} \colon C_n(M) \to S^{N}$ be the map which sends 
$(x_i) \in M^{\n} \subseteq (\R^{N+1})^{\n}$
to the unit vector in the direction of $x_j - x_i$. 

\item Let $I = [0,\infty]$, the one-point compactification of the nonnegative 
reals, and for $(i,j,k) \in C_3(\n)$ let $s_{ijk} \colon C_n(M) \to I = [0, \infty]$ 
be the map which sends $(x_i)$ to 
$\left(|x_i - x_j|/|x_i - x_k|\right)$. 

\item Let $A_n[M]$ be the product $M^{\und{n}} \times  (S^{N})^{C_2(\n)}
\times I^{C_3(\n)},$ and similarly
let $A_n\la M \ra = M^{\n} \times (S^N)^{C_2(\n)}$.

\item Let $C_n[M]$ be the closure of the image of $C_n(M)$ in $A_n[M]$ 
under $\alpha_n = \iota \times (\pi_{ij}) \times (s_{ijk})$.  

\item Let $C_n \la M \ra$ be the closure of
the image of
$C_n(M)$ in $A_n \la M \ra$ under $\beta_n = \iota \times (\pi_{ij})$.

\end{enumerate}
\end{definition}

We call $C_n[M]$ the canonical completion or compactification of $C_n(M)$.
It has alternately been called the Fulton-MacPherson or the Axelrod-Singer
completion.  The standard definition has been as a closure in a 
product of blow-ups of $M^S$ along their diagonals 
for $S \subset \n$ \cite{FM, AxelSing}, though
Gaiffi gives an elementary definition for Euclidean spaces similar to
ours in \cite{Gaif03}.  Because of \refC{cosimp} below, 
we call $C_n \la M \ra$ the simplicial
variant.  

\begin{theorem}\label{T:mainlist}
The completions $C_n[M]$ and $C_n \la M \ra$ have the following properties.
\begin{enumerate}
\item They are compact when $M$ is.
\item The inclusion $C_n(M) \to M^{\und{n}}$ factors through a  surjective map
$C_n[M] \to M^{\und{n}}$, or respectively $C_n\la M \ra \to M^{\und{n}}$.

\item The homeomorphism types of $C_n[M]$ and $C_n \la M \ra$
are independent of the ambient embedding in $\R^{N+1}$.
\item An embedding $f : M \hookrightarrow N$ functorially
induces a map   $ev_n(f) : C_n[M] \to C_n[N]$
(or by abuse $C_n \la M \ra \to C_n \la N \ra$), extending the induced maps on $C_n(M)$.
\item $C_n[M]$ is a manifold with corners with $C_n(M)$ as its interior. \label{4}
\item The maps $C_n(M) \hookrightarrow C_n[M] \overset{Q}{\to} C_n\la M \ra$, where
$Q$ is the restriction of the projection from $A_n[M] \to A_n \la M \ra$,
are homotopy equivalences. \label{heq}
\end{enumerate}
\end{theorem}

\begin{proof}
\begin{enumerate}
\item They are closed in $A_n [M]$ and $A_n \la M \ra$, respectively,
which are compact when $M$ is.
\item Restrict the projection of $A_n[M]$, respectively $A_n \la M \ra$, 
onto its factor $M^{\n}$.  We establish surjectivity only for compact $M$, where
it holds because the image
is closed (in fact compact) and contains an open dense subset (namely $C_{n}(M)$).
\item For $C_n[M]$, this is Theorem~4.7 of \cite{Sinh04}, which follows 
from careful analysis of $C_n[M]$ as a manifold with corners.  The statement for
$C_n \la M \ra$ follows from that for $C_n[M]$ and Theorem~5.8 of \cite{Sinh04}, which 
identifies $C_n \la M \ra$ as a pushout of $C_n[M]$ in terms which do not use
the embedding of $M$ in $\R^{N+1}$.
\item If $M \hookrightarrow N$, we can restrict the ambient embedding of $N$ in $\R^{N+1}$
to define the ambient embedding of $M$ in $\R^{N+1}$, so that $C_n[M]$ and $C_n\la M \ra$
are manifestly subspaces of $C_n [N]$ and $C_n \la N \ra$.  
\item This is Theorem~4.4 of \cite{Sinh04}, the culmination of a detailed local
analysis of $C_n[M]$.
\item That the inclusion $C_n(M) \to C_n[M]$ is a homotopy equivalence follows
from \ref{4} and the fact that a topological manifold with boundary is 
homotopy equivalent to its interior.  That the projection $C_n[M] \to C_n \la M \ra$
is an equivalence is 
Theorem~5.10 of \cite{Sinh04}, the main result of its Section~5.
\end{enumerate}
\end{proof}

We will use maps between these completed configuration spaces.
To map from $C_n[M]$ or $C_n \la M \ra$ it suffices to restrict a map from
$A_n[M]$ or respectively $A_n \la M \ra$.  
In defining maps to $C_n[M]$ or $C_n \la M \ra$ we may use explicit characterizations 
of them as subspaces of $A_n[M]$ or respectively $A_n \la M \ra$.  Such characterizations
are given in Theorems~4.1~and~5.14 of \cite{Sinh04}.

\subsection{Stratification and a category of trees}

\begin{definition}  
Define an $f$-tree to be a rooted, connected tree, with labeled
leaves, and with no bivalent  internal vertices.  Thus, an $f$-tree $T$ 
is a connected acyclic graph with a specified vertex $v_r$ called the root. 
The root may have any valence, but other vertices may not be bivalent. 
The univalent vertices, excluding the root if happens to be univalent, 
are called leaves, and  each leaf is labeled uniquely with an element of 
the appropriate $\n$.  
\end{definition} 

See Figure~\ref{F:Psi3} for examples of $f$-trees (noting that planar embedding
is not part of their structure).

\begin{definition} 
\begin{enumerate}
\item Given an $f$-tree $T$ and a set of non-leaf edges $E$ the 
contraction of $T$ by $E$ is the tree $T'$ obtained by, for 
each edge $e \in E$, identifying its initial vertex with its
terminal vertex and removing $e$ from the set of edges. 

\item Define $\Psi_{{n}}$ to be the category whose objects are $f$-trees with
$n$  leaves.  There is a unique morphism in $\Psi_{{n}}$ from $T$ to
$T'$ if $T'$ is  isomorphic to a contraction of $T$ along some set of edges. 
\end{enumerate}
\end{definition} 


\skiphome{
\begin{figure}\label{F:Psi3}
$$\includegraphics[width=6cm]{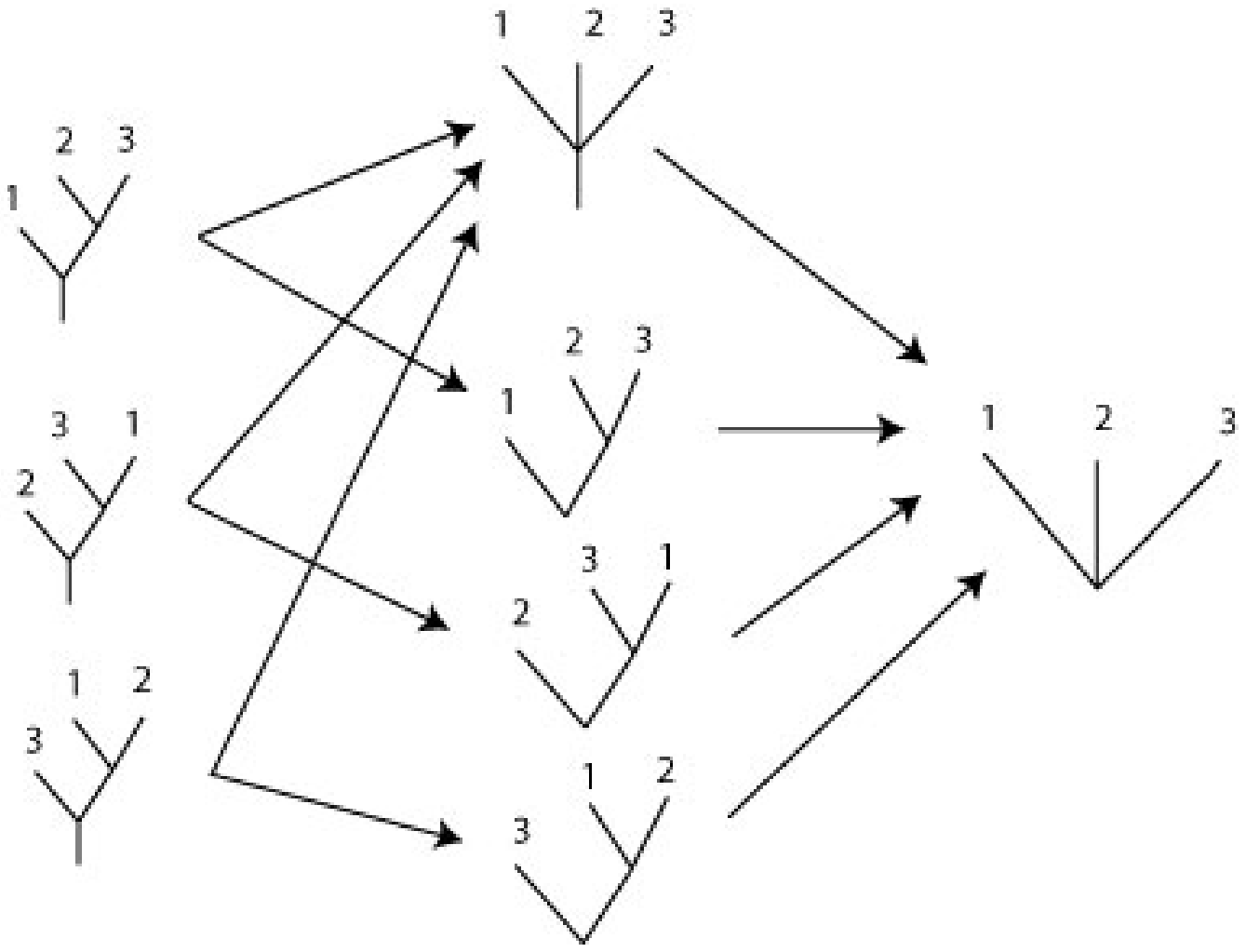}$$
\begin{center}
Figure~\ref{F:Psi3} - The category $\Psi_3$.
\end{center}
\end{figure}
}

This category of trees is essential in understanding the geometric 
structure of $C_n[M]$.   

\begin{definition}
\begin{enumerate} 
\item The root path of a leaf  is the unique 
path to the root. 
\item A leaf lies over a vertex if that vertex is in the leaf's root path.
\item Leaves $i$ and $j$ exclude leaf $k$ if there is some vertex $v$
such that $i$ and $j$ lie over $v$ but $k$ does not lie over $v$.
\item Given an $f$-tree $T$ with $n$ leaves, let $Ex(T) \subset \n^3$ be 
the subset of $(i,j,k)$ such that $i$ and $j$ exclude $k$ in $T$.
\end{enumerate}
\end{definition}

If  $Ex(T) = Ex(S)$ then either $T$ and $S$ are isomorphic or 
one is obtained from the other by adding an edge to the root vertex.

We now define a stratification of $C_n[M]$ through coordinates in $A_n[M]$
(recall our conventions about stratifications from Section~1).

\begin{proposition}
Let $x = (x_i) \times (u_{ij}) \times (d_{ijk})$ be a point in $C_n[M] \subset A_n[M]$.
There is a tree $T(x)$ such that $Ex(T(x))$ is the subset of $(i,j,k) \in \n^3$ such that
$d_{ijk} = 0$ .      We choose $T(x)$ to have a univalent root if and only if all of the $x_i$ are equal.  
\end{proposition}

See the comments after Definition~3.1 in \cite{Sinh04} for an indication of proof.

\begin{definition}
Given an $f$-tree $T$ let $C_T(M) \subset C_n[M]$ be the subspace of
$x$ such that $T(x) = T$, and let $C_{T}[M]$ be its closure.
\end{definition}

There is also a stratification of $C_{n}[M]$ arising from its structure as a manifold with
corners.  Namely, each stratum is a connected component
of the subspace of points modeled on $(\R^{\geq0})^{d_{1}} \times \R^{d_{2}}$
for fixed $d_{1}$ and $d_{2}$.

\begin{theorem}\label{T:strata}
The stratification of $C_n[M]$ by the $C_T(M)$ coincides with its stratification
as a manifold with corners.  The stratification poset is isomorphic to $\Psi_n$.  
Moreover, given an embedding $f$ of $M$ in $N$, the induced map
$ev_n(f) : C_n[M] \to C_n[N]$ preserves the stratification. 
\end{theorem}

This theorem summarizes some of the main 
results of \cite{Sinh04}: Theorem~3.4, which identifies
$\Psi_n$ as the poset associated to the $C_T(M)$ stratification;
Theorem~4.4, which establishes the manifold-with-corners structure on $C_n[M]$
and equates it with the tree stratification; 
and Theorem~4.8, which states that the evaluation map preserves strata.

We omit from this summary the explicit geometric description of $C_T(M)$
in terms of configurations in the tangent spaces of $M$, which
is the focus of much of \cite{Sinh04}.  We give the flavor of  how configurations degenerate to
give rise to configurations in the tangent bundle of $M$ in Figure~\ref{F:seq}.  
For readability we have omitted the labels
of points and leaves, necessary to work with ordered configurations.

\skiphome{
\begin{figure}\label{F:seq}
$$\includegraphics[width=12cm]{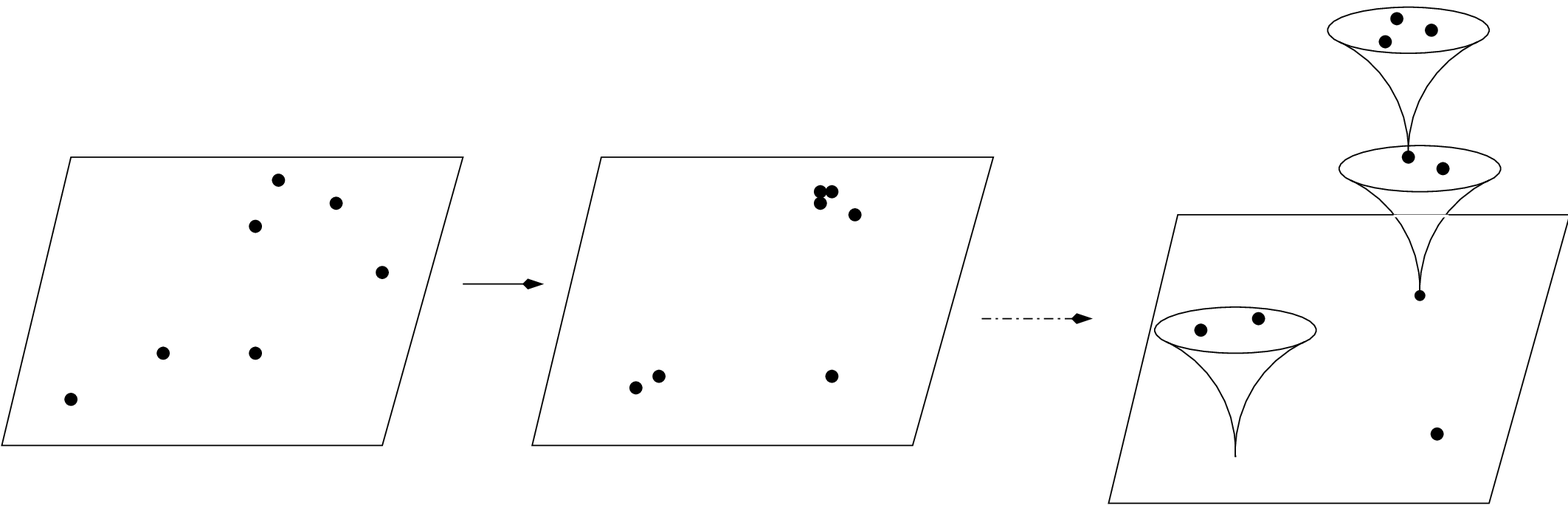}$$
\begin{center}
Figure~\ref{F:seq}. A sequence of points in $C_7(M)$  converging to 
a point in $C_{\includegraphics[width=0.4cm]{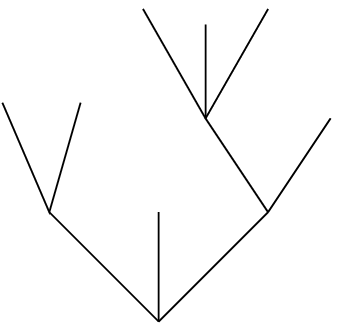}}(M)$.
\end{center}
\end{figure}
}


\subsection{Modifications to account for tangent vectors and boundary points}\label{S:tgtbdary}

We need tangent vectors at the constituent points of our configurations. 
We use $STM$ to denote the unit tangent bundle of $M$.

\begin{definition}
\begin{enumerate}
\item If $X(M)$ is a space equipped with a map to $M^{S}$, define
$X'(M)$ as the following pull-back 
$$
\begin{CD}
X'(M)  @>>> (STM)^{S} \\
@VVV @VVV\\
X(M) @>>> M^{S}.
\end{CD}
$$
\item If $f : X(M) \to Y(M)$ is a map over $M^{S}$, let $f' 
:  X'(M) \to Y'(M)$ be the induced map on pull-backs.
\end{enumerate}
\end{definition}

We next define a variant of these completions for manifolds with
selected tangent vectors in its boundary, say ${\bf v_0}, {\bf v_1} \in S(TM|_{\partial M})$
which sit over ${\bf y_0}$ and ${\bf y_1}$ in $\partial M$, as
needed in the definition of $\E$.  These points will be ``secretly added as first and last'' in all configurations, so we switch from indexing by $\n$ to indexing by $[n+1]$.

\begin{definition}
\begin{itemize}
\item Let $A_n[M, \partial]$ be the product $M^{[n+1]} \times  (S^{N})^{C_2([n+1])}
\times I^{C_3([n+1])},$ and similarly
let $A_n\la M, \partial \ra = M^{[n+1]} \times (S^N)^{C_2([n+1])}$.
\item Let $\iota^\partial : C_n(M - ({\bf y_{0} \cup y_{1}})) \to M^{[n+1]}$ send $(x_1, \ldots, x_n)$
to $({\bf y_0}, x_1, \ldots, x_n, {\bf y_1})$.  We use $\iota^\partial$ to define
$\pi_{ij}^\partial$ and $s_{ijk}^\partial$, following \refD{basiccomp}.
\item Let $C_n[M, \partial]$ be the closure of $C_n(M - ({\bf y_{0} \cup y_{1}}))$ in $A_n[M, \partial]$
under $\iota^\partial \times \pi_{ij}^\partial \times s_{ijk}^\partial$, and let
$C_n \la M, \partial \ra$ be defined similarly.
\item Let $(STM)^{n}_{\partial}$ be the subspace of $(STM)^{[n+1]}$ 
of collections whose first vector is  ${\bf v_0}$ and whose last is ${\bf v_1}$.
 Let $C_n'[M, \partial]$ be defined as the pullback
$$
\begin{CD}
C'_{n}[M, \partial] @>>> (STM)^{n}_{\partial} \\
@VVV @VVV\\
C_n[M, \partial]@>\iota^\partial>> M^{[n+1]},
\end{CD}
$$
and let $C_{n}' \la M, \partial \ra$ be defined similarly.
\end{itemize}
\end{definition} 

Through an obvious relabeling of coordinates, $C_{n}\la M, \partial \ra$ may be
viewed as a subspace of $C_{n+2} \la M \ra$, which we will need to do in 
Section~\ref{S:interp}.

\subsection{The associahedron}

When $M$ is an interval, the Fulton-MacPherson completion coincides
with the well-known associahedron.  
Note that $C_n(I)$ has one
component for every permutation of $n$ letters.   We focus on a single component.

\begin{definition}
Let $\wt{C}_n[I, \partial]$ denote the component of $C'_n[I, \partial]$
which is the closure of the component of $C'_n(I)$ for which the order of the points in the configuration agrees with the order they occur in the interval, and with all unit tangent vectors ``positive''.
\end{definition}

Passing to the ordered component of $C_n'[I, \partial]$ leads combinatorially to 
giving planar embeddings to the trees which label strata.

\begin{definition}
\begin{enumerate}
\item Call $S \subset \und{n}$ consecutive if
$i, j \in S$ and $i < k < j$ implies $k \in S$.  
\item  Let $\Psi^o_n$ denote the full sub-category of $\Psi_n$ whose objects
are $f$-trees such that the set of leaves over any vertex is consecutive
and such that the root vertex has valence greater than one.
\end{enumerate}
\end{definition}

Any element of $\Psi^o_n$ has an embedding in the upper half plane with 
the root at the origin, unique up to isotopy, in which labels of leaves coincide 
with their ordering given  by the clockwise orientation of the plane.   
We may then omit the labels of leaves from such an embedding.

Recall, for example from \cite{Stas70},  Stasheff's associahedron.  For
our purposes, a convenient definition is as the space of rooted half-planar
trees where edges have lengths and each leaf has a total distance of one from the
root.  We denote the
$n$-dimensional associahedron by $K_{n+2}$, using Stasheff's
terminology, as the relevant trees in this case have $n+2$ leaves.  
The following is Theorem~4.19 of \cite{Sinh04}, and is well-known
in different technical settings \cite{MaSt02}.

\begin{theorem}\label{T:stashdiff}
The completion $\wt{C}_n[I, \partial]$ and $K_{n+2}$
are diffeomorphic as manifolds with corners.
Moreover, their barycentric subdivisions are isomorphic as polyhedra to 
the realization (or order complex) of the poset $\Psi^o_{n+2}$.
\end{theorem}

\begin{proof}[Sketch of proof]
To understand $\wt{C}_n [I, \partial]$, we first use \refT{strata} to see that
its stratification poset is isomorphic to $\Psi^o_{n+2}$.  The remaining work is
mainly to show that the strata are all isomorphic to open Euclidean balls.

The relationship between $K_{n+2}$ and $\Psi^o_{n+2}$ is tautological
from our chosen definition.  For relations of this definition with others see
for example Proposition~1.22 of \cite{BoVo73}.
\end{proof}

A picture of $K_4$, the pentagon, as the realization of $\Psi^o_4$
is given in Figure~\ref{F2}.

In \cite{FeKo03} Feitchner and Kozlov view the posets $\Psi_{n+1}$  as a
combinatorial blow-ups of $\pnn$.   We relate these two posets in order 
to change indexing categories in the next section.

\begin{definition} \label{D:fn} 
\begin{enumerate}
\item We say two leaves  are root-joined if their root paths intersect only at the root vertex. 
\item We call the pair of leaves labeled
by $i$ and $i+1$ the $i$th adjacent pair. 
\item  Root-joined pairs of indices remain root-joined 
after applying a morphism on $\Psi^o_n$.  
Let $\F_n$ be the functor from $\Ston$ to $\pnn$ which sends a tree $T$ to 
the set $S$ where $i \in S$ if the $i$th adjacent pair of leaves in $T$
is root-joined.
\end{enumerate}
\end{definition} 

\skiphome{
\begin{figure}\label{F2}
\psfrag{F}{$\F_3$}
$$\includegraphics[width=12cm]{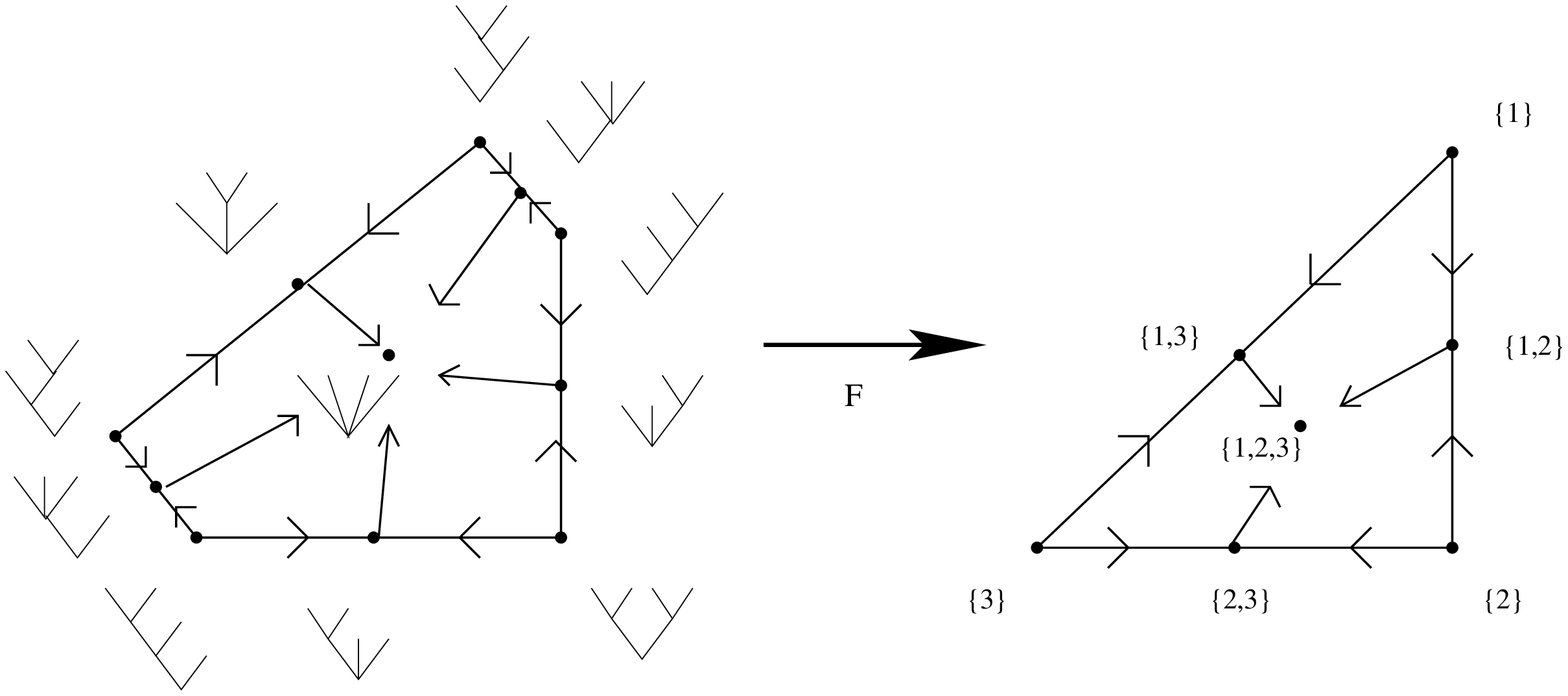}$$
\begin{center}
Figure~\ref{F2}. An illustration of the functor $\F_3$.
\end{center}
\end{figure}
}

See Figure~\ref{F2} for an illustration of $\F_3$.
Recall from \refP{holim} the notion of 
cofinality of a functor (between indexing categories).  

\begin{proposition}\label{P:cofin1} 
The functors $\F_n \colon \Ston \to \pnn$ are left cofinal. 
\end{proposition} 

\begin{proof} 
For any object $S$ of $\pnn$, the category $\F_n \downarrow S$ is a 
poset with a terminal object, namely the 
unique tree which maps to $S$ such that there are at most two
edges between any leaf and the root.  
\end{proof}

\subsection{Projection and diagonal maps}

Some of the most utilized maps between products of spaces are projection
and diagonal maps, which are a special case of maps defined when
one views $S \mapsto X^S$ as a contravariant functor from sets to spaces.

\begin{definition}\label{D:funct}
Given a map of sets $\sigma : R \to S$ let
$p^X_{\sigma}$, or just $p_{\sigma}$,
denote the map from $X^S$ to $X^R$ which sends 
$(x_i)_{i \in S}$ to $(x_{\sigma(j)})_{j \in R}$.
\end{definition}

In this section we define similar maps for completions of configuration spaces.
We start with $C'_n \la \R^{N+1} \ra$.  First observe that 
$A'_n\la \R^{N+1} \ra = (\R^{N+1} \times S^N)^{\und{n}} \times
(S^N)^{C_2(\und{n})}$ is canonically diffeomorphic to  $(\R^{N+1})^{\und{n}}
\times (S^N)^{\und{n}^2}$, letting $u_{ii}$ be the unit tangent vector
associated to the $i$th factor of $\R^{N+1}$.  Informally, we think of the tangent
vector $u_{ii}$ as ``the unit vector from $x_{i}$ to itself.''   We will use this data for
diagonal maps.

\begin{definition}\label{D:Fsigma}
Let $\sigma : \und{m}  \to \n$.  Define 
$A_\sigma : A'_n\la \R^{N+1} \ra \to A'_m \la \R^{N+1} \ra$ as 
$p^{\R^{N+1}}_\sigma \times p^{S^N}_{\sigma^2}$,
and let $F_\sigma$ be the restriction of $A_\sigma$ to $C_n\la M \ra$.
For the variants with boundary,
we say a map $\sigma : [m] \to [n]$ which
sends $0 \mapsto 0$ and $m \mapsto n$ is boundary-preserving.
For boundary-preserving $\sigma$,  define $A^{\partial}_{\sigma}$ and 
$F^{\partial}_{\sigma}$ in a similar manner as above.
\end{definition}

The following is essentially Proposition~6.6 of \cite{Sinh04}.

\begin{proposition}\label{P:Fsigma}
Given a boundary-preserving $\sigma : [m+1] \to [n+1]$
the induced map $F^{\partial}_\sigma$ sends 
$C'_n\la M, \partial \ra$ to $C'_m \la M, \partial \ra$
and commutes with $p_\sigma^{STM}$.
\end{proposition}

Let $\mathcal{N}$ be full subcategory of the category of sets whose objects are
the sets $\underline{n}$.  Let $\mathcal{N}_{\partial}$ be the category whose
objects are the sets  $[n]$ and whose morphisms are boundary-preserving
maps between $[m]$ and $[n]$.

\begin{corollary}
Sending $\underline{n}$ to $C_n' \la M \ra$ and $\sigma$ to $F_\sigma$ defines
a contravariant functor from $\mathcal{N}$ to spaces.
Similarly, sending $[n+1]$ to $C_n' \la M, \partial \ra$ and boundary-preserving
$\sigma$ to $F^{\partial}_\sigma$ defines
a contravariant functor from $\mathcal{N}_{\partial}$ to spaces.
\end{corollary}

\begin{proof}
We check that $F_{\sigma \circ \tau} = F_\tau \circ F_\sigma$ and 
$F^{\partial}_{\sigma \circ \tau} = F^{\partial}_\tau \circ F^{\partial}_\sigma$.  These follow from
checking the analogous facts for $p_\sigma$ and $p_{\sigma^2}$, which are
immediate.
\end{proof}

These constructions give projection maps when the underlying maps on coordiates
are  inclusions.  

\begin{proposition}\label{P:proj}
If $\sigma : [m+1] \to [n+1]$ is an inclusion, the square
$$
\begin{CD}
C'_{n}(M)  @>>> C'_{n} \la M, \partial \ra \\
@V\bar{p}_{\sigma}VV   @VF_{\sigma}VV \\
C'_{m}(M) @>>> C'_{m} \la M, \partial  \ra
\end{CD}
$$
commutes, where the horizontal maps are the canonical inclusions and 
$\bar{p}_{\sigma}$ denotes the restriction of $p_{\sigma}$ to configuration spaces.
\end{proposition}

There is a similar commuting square in the case without boundary.

We use these $F^{\partial}_\sigma$ as the structure maps in our cosimplicial model.
We set notation for cosimplicial spaces and translate
between the two standard ways to describe maps between simplices as follows.

\begin{definition}\label{D:basiccosimp}
\begin{enumerate}
\item Let $\bf{\Delta}$ have one object for each nonnegative $n$ and
whose morphisms are the non-decreasing ordered set morphisms between the 
sets $[n] = \{ 0, 1, \ldots, n \}$.  
A functor from $\bf{\Delta}$ to spaces is a cosimplicial space. 
\item Put coordinates on $\Delta^n$ by
$0 = t_0 \leq t_1  \leq \cdots \leq t_n \leq t_{n+1} = 1$, and label its vertices by
elements of $[n]$ according to the number of $t_i$ equal to one, not counting
$t_{n+1}$.
The  standard cosimplicial object $\Delta^\bullet$ sends $[n]$ to  $\Delta^n$ and
sends some $f: [m] \to [n]$ to the linear map which extends the map given by $f$ on vertices.  
\item Given an order-preserving
$\sigma: [n] \to [m]$, define $\sigma^* : [m+1] \to [n+1]$ by setting
$\sigma^*(j)$ to be the number of $i \in [n]$ such that $\sigma(i) >  m-j$.  Note that
$\sigma^{*}$ is boundary-preserving as well as order-preserving.
\end{enumerate}
\end{definition}
 
The significance of $\sigma^{*}$ is that 
the linear map extending $\sigma : [n] \to [m]$ on vertices
sends  $(t_i)_{i=0}^{n+1}  \in \Delta^n$ to $(t_{\sigma^*(j)})_{j=0}^{m+1}  \in \Delta^m$.

\begin{corollary}\label{C:cosimp}
The functor which sends $[n]$ to $C'_n \la M, \partial  \ra$ and 
$\sigma: [n] \to [m]$ to $F^{\partial}_{\sigma^{*}}$
defines a  cosimplicial space, denoted  $\C^\bullet \la M \ra$. 
\end{corollary}

When $M = I$ we have $\C^\bullet \la I \ra = \Delta^\bullet$, which is the reason
we call $C_n \la M \ra$ the simplicial variant completion of configuration spaces.
For explicit use later, it is helpful to name the coface maps 
of $\C^\bullet \la M \ra$.

\begin{definition}\label{D:cofaces}
Let $\delta^i : C'_n \la M, \partial \ra \to C'_{n+1}\la M, \partial \ra$ be the $i$th
coface map of $\C^\bullet \la M \ra$.  Explicitly, $\delta^i = F^{\partial}_{\tau_{i}}$ 
where $\tau_{i} :[n+1]\to[n]$ is
the order-preserving surjection in which $i$ and $i+1$
both map to $i \in [n]$.
\end{definition}

Informally we say  $\delta^i$ doubles the $i$th point in a configuration, thus
acting as a diagonal map.  

\medskip

For $C_n[M]$ projection maps work in a straightforward manner, but
we do not use them.  Diagonal maps are more involved, as they
will never satisfy the identities $\delta^i \circ \delta^i 
= \delta^{i+1} \circ \delta^i$;  see
the last comments of \cite{Sinh04} for an illustration.  However, the various
composites of diagonal maps do have canonical homotopies between them which 
are parameterized by associahedra.  We include these associahedra as part
of the basic definition of  diagonal maps.  We will start with composites of the 
diagonal maps we defined on $C_{n} \la M \ra$, as subspaces of $A_{n} \la M \ra$.
To go from $A_{n} \la M \ra$ to $A_{n} [ M ]$, we need to address
factors of $I^{C_{3}(\n)}$.  

Let $(e_{j \ell m}(a))$ denote the image of $a \in K_{k+1} = C_{k-1}[I, \partial]
\subset A_{k-1}[I, \partial]$ under the projection to $I^{C_3([k+1])}$.  This map
from $K_{k+1}$ to $I^{C_{3([k+1])}}$ is an embedding of the associahedron in the cube
which extends the map from $C_{k-1}(I)$ sending $(t_{i})$ to
$(e_{j \ell m} = \frac{t_\ell -t_j}{t_m - t_j})$.

\begin{definition}
\begin{enumerate}
\item Let  $[k]+i  = \{i, i+1, \ldots, i+k\}$, and 
let $\varphi_{i,k} : [n+k+1] \to [n+1]$ be the unique order-preserving 
surjection sending all of $[k] + i$ to $i$.  We use the notation $\varphi$ instead
of $\varphi_{i, k}$ when $i$ and $k$ are determined by context.

\item Define $\iota_{i, k} : I^{C_3([n+1])} \times K_{k+1} \to I^{C_3([n+k+1])}$
by  sending  $(d_{n o p}) \times a$ to
$(f_{j \ell m})$ with 
$$
f_{j \ell m} =
\begin{cases}
d_{\varphi(j) \varphi(\ell) \varphi(m)} & {\text{if at most one of}} \; j, \ell, m \in [k]+i \\
0 & {\text{if}} \; j,\ell \in [k]+i \; {\text{but}} \; m \notin [k]+i \\
1 & {\text{if}} \; \ell, m \in [k]+i \; {\text{but}} \; j \notin [k]+i \\ 
\infty & {\text{if}} \; j,m \in [k]+i \; {\text{but}} \; \ell \notin [k]+i \\
e_{j-i, \ell-i, m-i}(a) & {\text{if}} \; j, \ell, m \in [k]+i.
\end{cases}
$$
\item Let $D_{i,k} : A_n'[M, \partial] \times K_{k+1} \to A_{n+k}' [M, \partial]$ be the product of
$A_{\varphi_{i,k}} : A_n' \la M, \partial \ra \to A_{n+k}' \la M, \partial \ra$ with $\iota_{i,k}$.
Let $\delta^i(k)$ denote the restriction of $D_{i,k}$ to $C_n'[M] \times K_{k+1}$.
\end{enumerate}
\end{definition}

\begin{proposition}\label{P:deltai}
$\delta^i(k)$ sends  $C_n'[M, \partial] \times K_{k+1}$ to $C_{n+k}'[M, \partial]$.
\end{proposition}

The analogous result without modifications at boundary
points  is Proposition~6.11 in \cite{Sinh04}.
The idea of proof is to check that the image of
$\delta^i(k)$ is in $C_{n+k}'[M, \partial]$ using the conditions of Theorem~4.1 of 
\cite{Sinh04}, starting with the fact that points in $C_n'[M, \partial]$ satisfy
these conditions.  

Informally, $\delta^i(k)$ repeats the $i$th point in a configuration $k$ times.
Such repeating is allowed in the canonical completion $C'_n[M]$
as long as there is consistent data to distinguish the points
``infinitesimally''.  In this case these repeated points are
aligned in the direction of the unit tangent vector 
associated to the $i$th point in the configuration.
How they sit on that ``infinitesimal line'' is given by the coordinates of the
associahedron $K_{k+1}$ through its identification with
$C_{k-1}[I, \partial]$.

By specializing to $M=I$ and applying
\refT{stashdiff} we obtain
$\delta^i(k) : K_n \times K_{k+1} \to K_{n+k}$.  We use these maps to refine \refT{stashdiff} 
by explicitly identifying the strata of the associahedron and their inclusion maps. 

\begin{definition}\label{D:Astrat}
\begin{enumerate}
Let $K$ be the functor from $ \Psi^o = \bigcup_{n} \Psi^{o}_{n}$ to spaces
defined as follows.
\item On objects $K$ sends a tree $T$ to $K_T = \prod_{v} K_{|v|}$, where 
$v$ runs over all non-leaf
vertices of $T$ and $|v|$ is the number of outgoing edges of $v$.  
\item To the basic morphism 
which contracts the $i$th edge (in the planar ordering) over a vertex
$v$, say connecting $v$ to $v'$ which  under identification both become $w$,
$K$ associates the product of
$\delta^i(|v'| -1) : K_{|v|} \times K_{|v'|} \to K_{|w|}$ with
the identity maps on factors $K_s$ with $s \neq v, v', w$.
\item Extend $K$ to all morphisms in $\Psi^o$ by factoring into basic
morphisms.  The resulting maps are independent
of this factorization.
\end{enumerate}
\end{definition}

\begin{theorem}
The functor $K$ restricted to $\Psi^{o}_{n+2}$
realizes the stratification of $K_{n+2} = C_n[I, \partial]$ as a manifold with corners.
In particular,  $C_T[I, \partial] \cong K_T$.
\end{theorem}

The proof of this theorem follows that of Theorem~4.19 of \cite{Sinh04}.

Finally, we will need a compatibility result for diagonal
maps on the canonical and simplicial variant completions.  
Because we used the diagonal maps on $C_{n}' \la M, \partial \ra$
in our definition of those for $C_{n}'[M, \partial]$, we have the following.

\begin{proposition}\label{P:deltacomm}
The diagram
$$
\begin{CD}
C'_n[M, \partial] \times K_{k+1} @>>{\delta^i(k)}> C'_{n+k}[M, \partial] \\
@V{Q' \circ p_1}VV    @VV{Q'}V \\
C'_n\la M, \partial \ra @>>{(\delta^i)^{\circ k}}> C'_{n+k} \la M, \partial \ra
\end{CD}
$$
commutes, where  $Q'$ is the restriction of the quotient map from 
$A'[M]$ to $A'\la M \ra$ and $p_1$ is the projection onto the first factor.
\end{proposition}

\section{The mapping space model}\label{S:mapping}  

In order to define the mapping space model, we introduce sub-strata of
the standard stratification of $C'_{n}[M, \partial]$.

\begin{definition}\label{D:align}
Let  $x = (x_i, v_i) \times (u_{ij}) \times (s_{ijk})$ be a point in  
$C'_T[M, \partial]$, where
$x_i$ and $v_i$ denote a point in $M$ and a unit vector in $T_{x_i}M$ 
respectively.  Such a point is called aligned (with respect to $T$)
if for all $i,j$ which are not root-joined 
(which implies $x_i = x_j$), we have $v_i = v_j$
and $u_{ij}$ is the image of $v_i$ under the Jacobian of the embedding
of $M$ in $\R^{N+1}$.  We call the subspace of aligned points
of $C'_{T}(M, \partial)$ the 
aligned (sub-)stratum.  We denote it by $C_T^\alpha(M, \partial)$ and
its closure by $C_{T}^{\alpha}[M, \partial]$.
\end{definition}

\skiphome{
\begin{figure}\label{F:align}
$$\includegraphics[width=5cm]{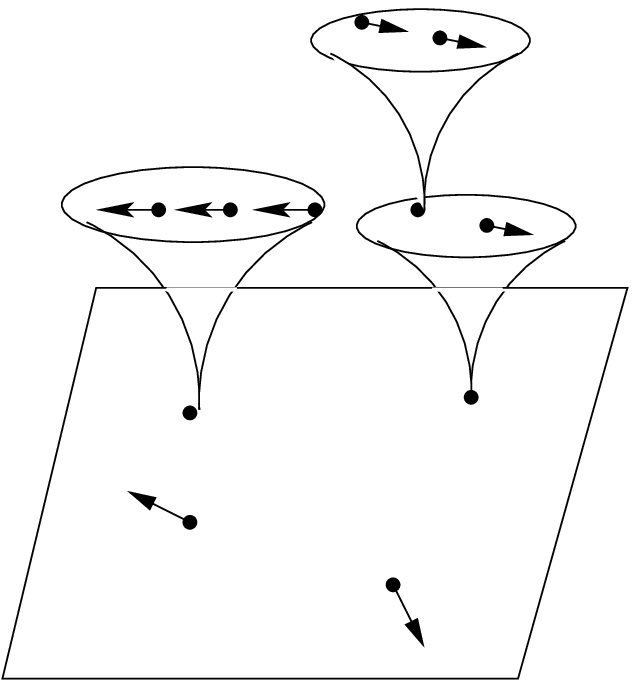}$$
\begin{center}
Figure~\ref{F:align} - A representation of an aligned point in $C_n'[M]$.
\end{center}
\end{figure}
}

We may build a tree inductively (non-uniquely)
as $T_{0} \subset T_{1} \subset \cdots \subset T_{\ell} = T$
where $T_{0}$ is the root and all of its outgoing edges and
$T_{j}$ is obtained from $T_{j-1}$ by adding  $m_{j}$ new edges which emanate 
from the leaf labeled $i_{j}$, which becomes the internal vertex $v_{j}$.

\begin{proposition}\label{P:idalign}
With  $T_{0} \subset T_{1} \subset \cdots \subset T_{\ell} = T$ as above, 
the aligned substratum $C_{T}^{\alpha}[M, \partial]$
is the image of $C'_{\# v_{r} - 2} [M, \partial] \times \prod_{v \in T, v \neq v_{r}} K_{\# v}$
under the composite 
$$\cdots \circ \left(\delta^{i_{j}} (m_{j} - 1) \times id \right)
         \circ \cdots \circ \left( \delta^{i_{1}} (m_{1} - 1)  \times id\right), $$
where $id$ denotes the identity map on factors of $K_{\#v}$ for $v$ other than $v_{j}$.
\end{proposition}

Proving this proposition is a straightforward matter of reconstructing
a point in $C'_{\# v_{r} - 2} [M, \partial] \times \prod_{v \in T, v \neq v_{r}} K_{\# v}$
from one in $C_{T}^{\alpha}$.   
Informally, the proposition  says that an aligned
configuration is determined by a configuration of fewer points to which diagonal or
``repeating'' maps are applied.  

By \refT{strata}, the evaluation map of a knot $\theta$, namely 
$ev_n(\theta)$ from $K_{n+2} =C_n[I,\partial]$ to
$C'_n[M, \partial]$ preserves the standard stratification.  Because a knot
is one-dimensional, configurations which degenerate
as points approach one another along a knot have
limits which lie in aligned substrata.

\begin{definition} 
A  map from $\widetilde{C}'_{n}[I, \partial]$ to $C'_n[M, \partial]$ which respects 
the sub-stratification by $C^\alpha_T[M, \partial]$ is called aligned.  
Let $AM_n(M)$ denote the space of 
aligned maps from $\widetilde{C}'_{n}[I, \partial] = K_{n+2}$ to $C'_n[M, \partial]$. 
\end{definition} 

We have that $AM_n(M)$ maps to $AM_{n-1}(M)$ by restricting an aligned map to a 
chosen principal face $K_T$ where $T$ has a single trivalent internal vertex. 
These restrictions are fibrations.  Let $AM_\infty(M)$ denote the inverse limit of 
the $AM_n(M)$. As noted above, the evaluation map $ev_n$ maps $\E$ to 
$AM_n(M)$. The $ev_n$ commute with restriction maps, so we let $ev_\infty$ from $\E$ to 
$AM_\infty(M)$ be their inverse limit.  

Recall from \refD{defpn} that $P_n \E$ 
denotes the $n$th polynomial approximation to the space of knots in $M$. 
Our first main result of this paper is the following. 

\begin{theorem}\label{T:mapping} 
$AM_n(M)$ is weakly homotopy equivalent to $P_n \E$ for all $n$ including
$n=\infty$.  Moreover, the evaluation map $ev_n$ coincides with the map 
$\alpha_n$ in the homotopy category. 
\end{theorem} 

From this theorem and Corollary~\ref{C:weq} due to Goodwillie and Klein,
we have the following.

\begin{theorem}\label{T:mapweq} 
Let the dimension of $M$ be greater than three.  The map $ev_\infty 
\colon \E \to AM_\infty(M)$ is a weak homotopy equivalence. 
\end{theorem} 

The rest of this section is devoted to proving 
\refT{mapping}.  We first identify the mapping space model as 
a homotopy limit and then find a zig-zag of equivalences between 
this homotopy limit and the homotopy limit defining $P_n \E$.  This
sequence first involves changing the shape of the homotopy limit
involved  and then interpolating between configuration spaces and 
embedding  spaces through a space which incorporates both. 

\subsection{Identifying $AM_n$ as a homotopy limit}

Recall the functor $K$ from \refD{Astrat}.  If there is a morphism between $T$ and $T'$, let 
$i_{T,T'}$ denote the corresponding inclusion of $K_T$ as a face of 
$K_{T'}$. Let $K_T^{nr}$ denote the product over non-root vertices
$\prod_{v \neq v_r} K_{|v|}$.  The
inclusion  $i_{T,T'}$ descends to a map, which we will by abuse give the same 
name, from $K_T^{nr}$ to $K_{T'}^{nr}$. 

\begin{definition} 
Define  $\D_n[M]$ to be the functor from $\Sto$ to spaces as follows:
\begin{itemize}
 \item A tree $T$ is sent to $C'_{|v_r|-2}[M, \partial] \times K_T^{nr}$.
 \item  A morphism $T \to T'$ which is contraction of 
a non-root edge is sent to the product of the identity
map on ${C'_{|v_r|-2}[M, \partial]}$  with $i_{T,T'}$.
\item A morphism $T \to T'$ which is the  contraction of the $i$th root edge 
$e$ of $T$ is sent to  
$$\delta^i(|v_t| - 1) \times id \colon\left( \left( C'_{|v_r|-2}[M, \partial]
\times K_{|v_t|} \right) \times \prod_{\substack{v \in T \\ v \neq v_r, v_t}} K_{|v|} \right)
\to  \left( C'_{|v_r'|-2}[M,  \partial]  \times  \prod_{\substack{v \in T' \\ v \neq v_r}} K_{|v|} \right),$$ 
where $v_t$ is the terminal (that is, non-root) vertex of $e$ and $v_r'$ is the root
vertex of $T'$.
\end{itemize}
\end{definition} 

When $M = I$, $\D_n[M]$ is the functor $K$ of \refD{Astrat} restricted to $\Psi^{o}_{n+2}$.

\begin{lemma}\label{l-1}
The space of aligned maps $AM_{n}(M)$ is homeomorphic to the 
homotopy limit of $\D_n[M]$. 
\end{lemma} 

This lemma follows from  
Proposition~\ref{P:inclimit} which identifies homotopy limits such
as $\D_n[M]$ with spaces of stratum-preserving maps,  Theorem~\ref{T:stashdiff}
which identifies $C_{n}[I, \partial]$ with the realization of $\Psi^o_{n+2}$,
and Proposition~\ref{P:idalign} which identifies 
the images of composites of $\delta^i(k)$ with aligned substrata.

\subsection{Changing from the canonical to the simplicial completion}

We begin to  interpolate between $\D_n[M]$, whose homotopy
limit is the mapping space model,  and $\Em_n(M)$, whose homotopy
limit is the degree~$n$ embedding calculus approximation to the
space of knots in $M$.  

\begin{definition} 
\begin{enumerate}
\item Let $\wt{\D}_n \la M \ra$ be the functor from $\Sto$ to spaces which 
sends $T$ to $C'_{|v_r|-2}\la M, \partial \ra$, sends the
contraction of the $i$th root edge to $(\delta^i)^{\circ (|v_t| - 1)}$ where $v_t$
is the terminal (non-root) vertex of the contracted edge, and sends the 
contraction of a non-root edge to the identity map.
\item Let $Q^{\D}_n : \D_n [M] \to \wt{\D}_n \la M \ra$ send 
$C'_{|v_r|-2}[M, \partial] \times K_T^{nr}$
to $C'_{|v_r|-2}\la M, \partial \ra$ through projection onto
$C'_{|v_r|-2}[M, \partial]$ composed with $Q'$,  where $Q'$ is as
in \refP{deltacomm}.
\end{enumerate}
\end{definition} 

By \refP{deltacomm}, $Q^{\D}_n$ is a natural transformation.
By part~\ref{heq} of \refT{mainlist}, $Q : C_n[M] \to C_n \la M \ra$ is a homotopy
equivalence, as is $Q'$ since it is induced by $Q$ on pull-backs.  Because
$Q^{\D}_n$ is an objectwise weak equivalence, we deduce the
following.

\begin{lemma}\label{l0} 
The map from the homotopy limit of $\D_n[M]$ to the homotopy limit of 
$\wt{\D}_n \la M \ra$ defined through $Q^{\D}_n$  is a weak 
equivalence. 
\end{lemma} 

\subsection{Changing the indexing category from trees to subsets}

Many of the maps in the definition of $\wt{\D}_n \la M \ra$ are the identity,
so they  may be eliminated, replacing
the category $\Sto$ by $\pnno$.  

\begin{definition}\label{D:dn} 
Let $\D_n \la M \ra$ be the functor from $\pnno$ to spaces which sends $S$ 
to $C'_{\#S-1}\la M, \partial\ra$  and sends the inclusion $S \subseteq S'$ 
with $S' = S \cup k$ to $\delta^i$, where $i$ is the number of elements 
of $S$ less than $k$.
\end{definition} 

By construction, we have the following. 

\begin{lemma} 
$\wt{\D}_n \la M \ra $ is the composite of $\D_n \la M \ra$ with the 
functor $\F_n$ of Definition~\ref{D:fn}. 
\end{lemma} 

The following lemma, which is the next link in our chain of equivalences,
is immediate from Proposition~\ref{P:cofin1},
which says that $\F_{n}$ are left cofinal, and 
the fact that left cofinal functors induce equivalences on homotopy limits. 

\begin{lemma}\label{l1}
The homotopy limit of ${\D}_n \la M \ra$ is weakly homotopy equivalent 
to the homotopy limit of $\wt{\D}_n \la M \ra$. 
\end{lemma} 

\subsection{Interpolating between configurations and punctured knots}\label{S:interp}

To interpolate between $\D_n \la M \ra$ and $\Em_n(M)$, we incorporate 
both embeddings and configurations in one space.   

\begin{definition} 
Given a metric space $X$ define $\Ha(X)$ to be the space whose points are 
compact subspaces of $X$ and with a metric defined as follows.  Let $A$ 
and $B$ be compact subspaces of $X$ and let $x$ be a point in $X$. 
Define $d(x,A)$ to be $\inf_{a \in A} d(x,a)$.  Define 
the Hausdorff metric $d(A,B)$ to 
be the greater of $\sup_{b \in B} d(b, A)$ and $\sup_{a \in A} d(a, 
B)$. 
\end{definition} 
 
Because $STM$ is metrizable, as is of course $S^{N}$, so are $C'_n \la M \ra$ and
$C'_n \la M, \partial \ra$.  Recall Definition~\ref{D:ej}, where $\{J_i\}$, for $i \in \no$ is 
a fixed collection of sub-intervals of $I$ and $E_S(M)$ is the space of constant speed
embeddings of the punctured interval $(I - \bigcup_{s \in S} 
J_s)$ in $M$. Recall as well from Section~\ref{S:tgtbdary} that $C_{n-1}\la M, \partial \ra$
may be viewed through a simple relabeling of coordinates as a subspace of
$C_{n+1}\la M \ra$.

\begin{definition} 
\begin{enumerate}
\item Fix some $S \subseteq \no$.
Let $\{ I_k \}$ be the set of 
connected components of $I - \bigcup_{s \in S} J_s$, and given an embedding $f$
of $(I - \bigcup_{s \in S} J_s)$ let $f_k$ be the restriction of $f$ to 
$I_k$.  

\item For $f \in E_S(M)$ let $ev_S(f)$ be the evaluation map from  
$\prod_{k=0}^{\#S} I_k$ to $C'_{\#S + 1} (M)$.

\item Define $E_S\la M, \partial \ra$ to be the union of  $E_S(M)$ and 
$C'_{\#S -1}\la M, \partial \ra$ as a subspace of 
$\Ha\left(C'_{\#S + 1} \la M \ra\right)$, where 
$C'_{\#S -1}\la M, \partial \ra \subset \Ha\left(C'_{\#S + 1} \la M 
\ra \right)$ is a subspace of one-point subsets (appropriately reindexed)
and some $f \in E_S(M)$ 
defines a point in $\Ha\left(C'_{\#S + 1} \la M \ra \right)$ by the image of $ev_S(f)$. 
\end{enumerate}
\end{definition}   

Because points in $E_S(M)$ are determined by their image and
$E_S(M)$ is endowed with the $C^1$-topology, it indeed embeds as a
subspace of  $\Ha(C'_{\#S + 1} \la M \ra)$, 
which also accounts for the first derivative in its topology.
We think of $E_S\la M\ra$ as a space of embeddings which may 
be degenerate by having all of the embeddings of components ``shrink
until they become tangent vectors''.  

\begin{proposition}\label{P:incl} 
The inclusions of $C'_{\#S-1}\la M, \partial \ra$ and  $E_S(M)$ 
in $E_S\la M\ra$ are weak homotopy equivalences.
\end{proposition}

We first define a retraction from $E_S\la M \ra $ onto  $C'_{\# S-1}\la 
M, \partial \ra$.  

\begin{definition} 
Let $m \in \prod I_k$ have first coordinate $m_0 = 0$, 
last coordinate $m_{\#S + 1} =1$, and other coordinates given by
defining $m_k$ to be the mid-point of  $I_k$. 
Let $\epsilon \colon E_S\la M\ra \to C'_{\# S-1}\la M, \partial \ra$  
be the identity map on $C'_{\# S - 1}\la M, \partial \ra$ and on 
$E_S(M)$ be defined by sending $f$ to
$(ev_S(f))(m)$.
\end{definition} 

\begin{proof}[Proof of Proposition~\ref{P:incl}] 
To show the inclusion of $C'_{\#S-1}\la M, \partial \ra$ is an equivalence
we show that $\epsilon$ is a deformation retraction.  Define the required
homotopy  on 
$E_S(M)$ by letting $s$ be the homotopy variable and setting  $f_k(t)_s = 
f((1-s)t + s \cdot m_k)$.   To show that the inclusion of $E_S(M)$ is
an equivalence, we show that the restriction of $\epsilon$ to $E_S(M)$
mapping to $C'_{\#S -1}(M, \partial)$ is an equivalence and then apply
part~\ref{heq} of \refT{mainlist} to get the equivalence with 
$C'_{\#S - 1} \la M, \partial \ra$.

This restriction of $\epsilon$  is a fibration by the isotopy
extension theorem.  We show that the fiber of this map, namely the space 
of embeddings of $I - \bigcup_{s \in S} J_s$ with prescribed  tangent vectors 
at the $m_k$, has trivial homotopy groups.  Suppose we 
have a family of such embeddings parameterized by a sphere.  
First we may apply a reparameterizing homotopy 
$G(t,s) = f((1-s)t + s m_k)$ for $s \in [0,a]$ for some $a$ so that 
the image of each component lies in a fixed Euclidean chart in $M$ about
the image of $m_k$.  By compactness there is an $a$ which works 
for the entire family.   

We choose coordinates in each chart around these points so that 
$f_k(t) = (f_{k,1}(t), f_{k,2}(t), \ldots)$ with each 
$f_{k,j}(m_k) = 0$ and $f_k'(m_k) = (b, 0, 0,
\ldots)$ for some $b>0$.  For each component consider the
``projection'' homotopy with $s \in [0,1]$ defined component-wise by 
$$H_{k}(t)_s = (f_{k,1}(t), s \cdot f_{k, 2}(t), s \cdot f_{k, 
3}(t), \dots).$$   
This homotopy is not necessarily a valid homotopy through embeddings, 
but it will always be on some 
neighborhood of $m_k$ since the derivative there is bounded away 
from zero throughout the homotopy.  By compactness, 
there is some non-zero $c$ such that this homotopy is an isotopy
on a neighborhood $N$ of $m_k$ of length $c$ for all points in the 
parameter sphere.  The composite of $G$, a 
second reparameterizing homotopy which changes the image of each interval 
so as to be the image of $N$, the projection homotopy $H$ and a rescaling
homotopy on the first coordinates of the fixed charts defines a 
homotopy between the given sphere of embeddings and 
a constant family.  Thus the fiber of $\epsilon|_{E_{S}(M)}$ is weakly contractible
as claimed.
\end{proof} 

The spaces $E_S\la M\ra$ are a suitable interpolation between  
$E_S(M)$ and $C'_{\# S-1}\la M, \partial \ra$.   We now define a 
suitable diagram interpolating between $\D_n\la M\ra$ and $\Em_n(M)$. 

\begin{proposition}\label{P:dbl} 
Let $S, S' \in \pnno$ with $S' = S \cup k$.  There is a map 
$\rho_{S,S'} \colon E_S\la M\ra \to E_{S'}\la M\ra$ whose restriction to  
$E_S(M)$ is the map to $E_{S'}(M)$ defined by restriction of embeddings 
and whose restriction to $C'_{\# S -1}\la M, \partial \ra$ is the map 
$\delta^i$ where $i$ is the number of elements 
of $S$ less than $k$. 
\end{proposition} 

\begin{proof} 
We check that the function so defined is continuous.  By definition
it is continuous when restricted to either $E_S(M)$ or $C'_{\# S -1}\la M, 
\partial \ra$.  Because $E_S(M)$ is open in $E_S\la M\ra$, it 
suffices to check continuity on  $C'_{\# S -1}\la M, \partial\ra$.  
Since $\delta^i$ is continuous on $C'_{\# S -1}\la M, \partial\ra$ itself,
it suffices to show that for every
$\varepsilon$ there is a $\gamma$ such that if distance between  $\theta \in 
E_S(M)$ and 
${\bf x} \in C'_{\# S -1}\la M, \partial\ra$ is less than $\gamma$ 
then the distance between 
$\omega = \rho_{S,S'}(\theta)$ and ${\bf y} = \delta^i({\bf x})$ is less than $\varepsilon$. 

Recall that  $C'_{\# S -1}\la M, \partial\ra$ is a 
subspace of $(STM)^{[\#S]} \times (S^{N})^{C_{2}([\#S])}$, 
whose points we label ${\bf x} = (v_j) \times (u_{jk})$.   Similarly let ${\bf y} =
(w_j) \times (z_{jk})$, where because ${\bf y} = \delta^i({\bf x})$ we have 
$w_j = v_{\tau_i(j)}$ and $z_{jk} = u_{\tau_i(j) \tau_i(k)}$, where
$\tau$ is as in \refD{cofaces}.
Let $\{I_k'\}$ be the components of $I - \{\bigcup_{s \in S'} J_s\}$.
Finally let $F$ be the ambient embedding of $M$ in $\R^{N+1}$ used to define $C_n \la M \ra$.  

A bound  on the distance between ${\bf x}$ and $\theta$ 
in $\Ha(C'_{\#S + 1} \la M \ra)$  is equivalent to a bound for each 
$k$ for the distance between $v_k$ and $\theta'(t)$ for all $t \in I_k$, 
as well as a bound on the distance between $u_{jk}$ and the unit vector 
in the direction of $\left(F \circ \theta(t) - F \circ \theta(s)\right)$ for $t \in
I_j$ and $s \in I_k$.  Such bounds give rise to the same  
bounds on the distance between $w_j$  and 
$\theta'(t)$ for all $t \in I_j$, as well as 
better bounds between the $z_{jk}$ and $\left(F \circ \theta(t) - F \circ \theta(s)\right)$
except when $j,k = i, i+1$.  In this last case, we may choose our bound on the 
distances between $v_i$ and $\theta'(t)$ for all $t \in I_i$ so that 
by the triangle inequality and elementary calculus the 
distance between $F'(w_i)$ and $(F\circ \theta(t_i) - F \circ
\theta(t_{i+1}))$ for all $t_i \in I_i, t_{i+1} \in I_{i+1}$ is 
arbitrarily small. 
\end{proof} 

\begin{definition} 
Let $\Em_n \la M\ra$ be the functor from $\pnno$ to spaces which sends $S$ 
to $E_S\la M\ra$ and sends the inclusion $S \subseteq S'$ where $S' = S 
\cup k$ to $\rho_{S,S'}$. 
\end{definition} 

We have constructed the maps $\rho_{S,S'}$ so that 
both $\D_n \la M \ra$ and $\Em_n(M)$ map to $\Em_n\la M \ra$ through 
inclusions entry-wise.  By Proposition~\ref{P:incl} these inclusions 
are weak equivalences, so we may deduce the following. 

\begin{lemma}\label{l3} 
The homotopy limits of $\D_n\la M\ra$  and $\Em_n(M)$ are weakly 
equivalent to the  homotopy limit of $\Em_n\la M\ra$. 
\end{lemma} 

\subsection{Assembling the equivalances}
We may now put together the proof of the main theorem of this section. 

\begin{proof}[Proof of Theorem~\ref{T:mapping}] 
 In Lemma~\ref{l-1} we showed that $AM_n(M)$ is 
homeomorphic to the homotopy limit of a diagram $\D_{n}[M]$.  
Recall that $P_n \E$ is the homotopy limit of the diagram $\Em_{n}(M)$.
The weak equivalence $AM_n(M) \simeq P_n \E$ then follows from the 
zig-zag of equivalences given by Lemmas~ \ref{l0}, \ref{l1}, 
and \ref{l3}, as assembled below.

\begin{equation} 
\D_{n}[M]   \underset{\ref{l0}}{\rightarrow} \wt{\D_{n}}\la M \ra 
\underset{\ref{l1}}{\leftarrow} \D_{n} \la M \ra  
\underset{\ref{l3}}{\rightarrow} \Em_{n} \la M \ra  
\underset{\ref{l3}}{\leftarrow} \Em_{n}(M). 
\end{equation}  

It remains to show that the evaluation map $ev_n : \E \to AM_n(M)$ coincides in the  
homotopy category with the map $\alpha_n$.  Clearly $ev_n$ coincides with 
other evaluation maps (which by abuse we also call $ev_n$) 
in the equivalences of Lemmas~\ref{l0} and \ref{l1}. 
We focus on the equivalences of Lemma~\ref{l3} and show that the 
following diagram commutes up to homotopy:
$$
\begin{CD}
{\rm Emb}(I, M) @>\alpha_{n}>> {\rm holim}\; \Em_{n}(M) \\
@V{ev_{n}}VV @VVV \\
{\rm holim}\;  \D_{n} \la M \ra  @>>> {\rm holim}\; \Em_{n} \la M \ra.
\end{CD}
$$

Recall that $\text{holim} \; \Em_n \la M\ra$ is a subspace of  
$(f_{S}) \in \prod_{S \in \pnno} 
 \text{Maps} (\Delta^{\#S -1}, E_S\la M \ra)$ consisting of products
 of maps compatible under restriction. 
The map $\alpha_n$ sends $\E$ to  $\text{holim} \; \Em_n \la M\ra$ 
as the subspace in which each $f_S$ is 
constant as a function on $\Delta^{\#S -1}$,  with image given by the 
restriction from $\E$ to $E_S(M)$. On the other hand, $ev_n$ maps $\E$
to the subspace $\text{holim} \D_n \la M \ra$ by evaluation.  
We will define a homotopy between these by 
``shrinking towards the evaluation points''.

We focus on the factor of  $\text{holim} \; \Em_n \la M\ra$ labeled by $S = \no$ itself.
Let $\tau = (t_i)$ be a point in  $\Delta^{n}$, 
and let $\theta \in \E$. 
Let $\rho(J, t, s)$ be the interval which linearly interpolates, with 
parameter $s$, between the interval $J$ and the degenerate interval 
$[t,t]$.  If  $J = [a,b]$ then $\rho(J, t, s)$ is the interval
$[(1-s) a + s t, (1-s)b + s t]$.  
Using adjointness, we define our homotopy through a map
$$ h_{\no} : \E \times \Delta^n \times I \to E_{\no} \la M \ra.$$ 
For $s<1$ set $h_{\no}(\theta, \tau, s)$  
to be the embedding of $\bigcup_k I_k$ which is the composite of 
the linear  isomorphism between $\bigcup_k I_k$ and 
$\bigcup_k \rho(I_k, t_k, s)$  and the restriction of
$\theta$ to $\bigcup_k \rho(I_k, t_k, s)$.  Define 
$h_{\no}(\theta, \tau, 1)$ to be $(ev_n(\theta))(\tau)$.  With the definitions
as given it is straightforward 
to check that $h_{\no}$ is  continuous.
For $S \subset \no$ we may lift $h_{\no}$ (through
piecewise linear homeomorphisms of intervals composed with $\theta$) to an
$h_{S}: \E \times \Delta^{\#S - 1} \times I \to E_{S} \la M \ra$, yielding 
a compatible family of homotopies as required.
\end{proof} 

Goodwillie and Klein's Theorem~\ref{T:good}, upon which we build, can be 
proved for knots in manifolds of dimension five or greater by 
dimension-counting arguments applied to the punctured embedding
spaces defining $\Em_{n}(M)$  (sharper versions of this theorem require 
surgery theory and the results of Goodwillie's thesis \cite{Go83}).  
It would be interesting to find a proof that the inclusion of the
knot space in the mapping space model is highly connected through
more direct arguments than the chain of equivalences just given,
perhaps through dimension-counting arguments.
Such an approach might
be useful for applications to classical knots, where Theorem~\ref{T:good} is not known.

\section{The cosimplicial model}\label{S:cosimp} 

We now produce cosimplicial models of  knot spaces.  We take as our 
starting point the model defined by $\D_n \la M\ra$, showing that 
it is pulled back from a cosimplicial diagram.  
Recall that the totalization of a cosimplicial space 
$\X^\bullet$ is the space of natural transformations from $\Delta^\bullet$ to $\X^\bullet$.
The following is ultimately the main theorem of this section. 

\begin{theorem}\label{T:cosimodel}
The space $\E$ is weakly equivalent to the totalization of a fibrant
replacement of the cosimplicial space $\C^\bullet \la M \ra$, defined
in \refC{cosimp}.
\end{theorem}

Recall the definition of the cosimplicial category $\bf{\Delta}$ from \refD{basiccosimp}.    

\begin{definition}
Let $\bf{\Delta}_n$ be the full subcategory
of $\bf{\Delta}$ whose objects are the sets $[i]$ for $0 \leq i \leq n$. 
Given a cosimplicial space 
$\X^\bullet$ let $i_n \X^\bullet$ be the restriction of $\X^\bullet$ to 
$\bf{\Delta}_n$.   
\end{definition}


The fact that the nerve of $\pnno$ is isomorphic to the barycentric 
subdivision of an $n$-simplex reflects the existence of a canonical 
functor from $\pnno$ to $\bf{\Delta}_{n}$. 

\begin{definition}\label{D:gn} 
Let $\G_n \colon \pnno \to \bf{\Delta}_{n}$ be the functor which sends a subset 
$S$ to $[\#S -1]$ and which
sends and inclusion $S \subseteq S'$ to the composite $[\#S - 1]
\cong S \subset S' \cong [\#S' - 1]$.  
\end{definition}

The maps in the cubical diagram $i_{n} \X^{\bullet} \circ \G_{n}$ are
the coface maps of $\X^{\bullet}$.

The proof of the following proposition is
an unraveling of Definitions~\ref{D:dn} and 
\ref{D:gn} and the definition of $\C^\bullet \la M \ra$ in \refC{cosimp}. 

\begin{proposition}\label{P:comp} 
$\D_n \la M \ra = i_{n} \C^\bullet \la M \ra \circ \G_n$. 
\end{proposition} 

The next step in proving Theorem~\ref{T:cosimodel} is to establish a 
general theorem about cosimplicial spaces and subcubical diagrams 
which is well-known to some but
for which we have not found a reference.  

\begin{theorem}\label{T:cusimp} 
Let $\X^\bullet$ be a cosimplicial space.  The homotopy limit of 
$i_{n} \X^\bullet \circ\G_n$ is weakly equivalent to the $n$th totalization of a
fibrant replacement of $\X^\bullet$.
\end{theorem}

Before proving this theorem in general, it is enlightening to establish
its first nontrivial case more explicitly.  Consider the homotopy limit
$$ H = {\text{holim}} \; X_0 \overset{d_0}{\to}
X_1 \overset{d_1}{\leftarrow} X_0,$$  
where $X_0$ and $X_1$ are entries of a fibrant cosimplicial space
$\X^\bullet$ with structure maps $d_0, d_1$ and $s_0$.
By definition, $d_0$ $d_1$ are sections of $s_0$, which is a
fibration.  We claim that this homotopy limit is weakly equivalent to
$\text{Tot}_1(\X^\bullet)$.  The 
homotopy limit $H$ naturally fibers over ${X_0}^2$ with fiber $\Omega
X_1$, the based loop space of $X_1$.  On the other hand, the first total
space fibers over $X_0$ with fiber equal to $\Omega( {\text{fiber}} \;
\pi)$, so the equivalence is not a triviality.  

Considering the diagram
$$
\begin{CD}
X_0 @>d_0>> X_1 @<d_1<<  X_0 \\
@VVidV        @Vs_{0}VV         @VidVV \\
X_0 @>id>> X_0 @<id<<  X_0
\end{CD}
$$
we see $H$ also fibers over $H_0$, the homotopy limit
of $X_0 \overset{id}{\to} X_0 \overset{id}{\leftarrow}
X_0$ which is homotopy equivalent to
$X_0$ through a deformation retraction onto the constant paths. The
fiber of this map over a constant path is homotopy equivalent to $\Omega (
{\text{fiber}} \; \pi)$.  In fact if we lift the homotopy equivalence of
$H_0$ with $X_0$ defined by shrinking a path to a constant path, we get a
homotopy equivalence of $H$ with a subspace of $H$ which is homeomorphic
to  ${\rm Tot}_1 \X^\bullet$.  

\medskip

In the general case \refT{cusimp} is immediate from two theorems, the first
of which is due to Bousfield and Kan \cite{BK}.

\begin{theorem}\label{T:BK}
The homotopy limit of $i_n \X^\bullet$ is weakly
equivalent to the $n$th totalization of a fibrant replacement of
$\X^\bullet$.
\end{theorem}

\begin{theorem}\label{T:Gcofin}
The functor $\G_n$ is left cofinal.
\end{theorem}

These two theorems along with Proposition~\ref{P:comp} give a chain of
equivalences which establishes Theorem~\ref{T:cosimodel}:
$$
{\text{holim}} \D_n \la M \ra  \underset{\ref{P:comp}}{\cong}  {\text{holim}}
 i_n
\C^\bullet \la M \ra \circ \G_{n }\underset{\ref{T:Gcofin}}{\leftarrow}
{\text{holim}} i_n
\C^\bullet \la M \ra  \underset{\ref{T:BK}}{\simeq} {\widetilde{\text{Tot}}}_n(\C^\bullet
\la M \ra ),
$$
where ${\widetilde{\text{Tot}}}_{n} \X^{\bullet}$ denotes the $n$th totalization of a 
fibrant replacement of $\X^{\bullet}$.
We proceed to prove Theorem~\ref{T:Gcofin}.  We say a simplicial complex is $n$-dimensional
if each simplex is a faces of some $n$-simplex.

\begin{lemma}\label{L:l1}
The realization $|\G_{n} \downarrow [d]|$ is isomorphic to the 
barycentric subdivision of the $n$-dimensional 
simplicial complex $Y_{n\to d}$ whose $i$-simplices 
are indexed by pairs $(S, f)$ where $S \subseteq [n]$ is of 
cardinality $i+1$ and 
$f$ is an order preserving map from $S$ to $[d]$, and whose face 
structure is defined by restriction of the maps $f$. 
\end{lemma} 

\begin{figure}\label{F:Ydn}
$$\includegraphics[width=7.2cm]{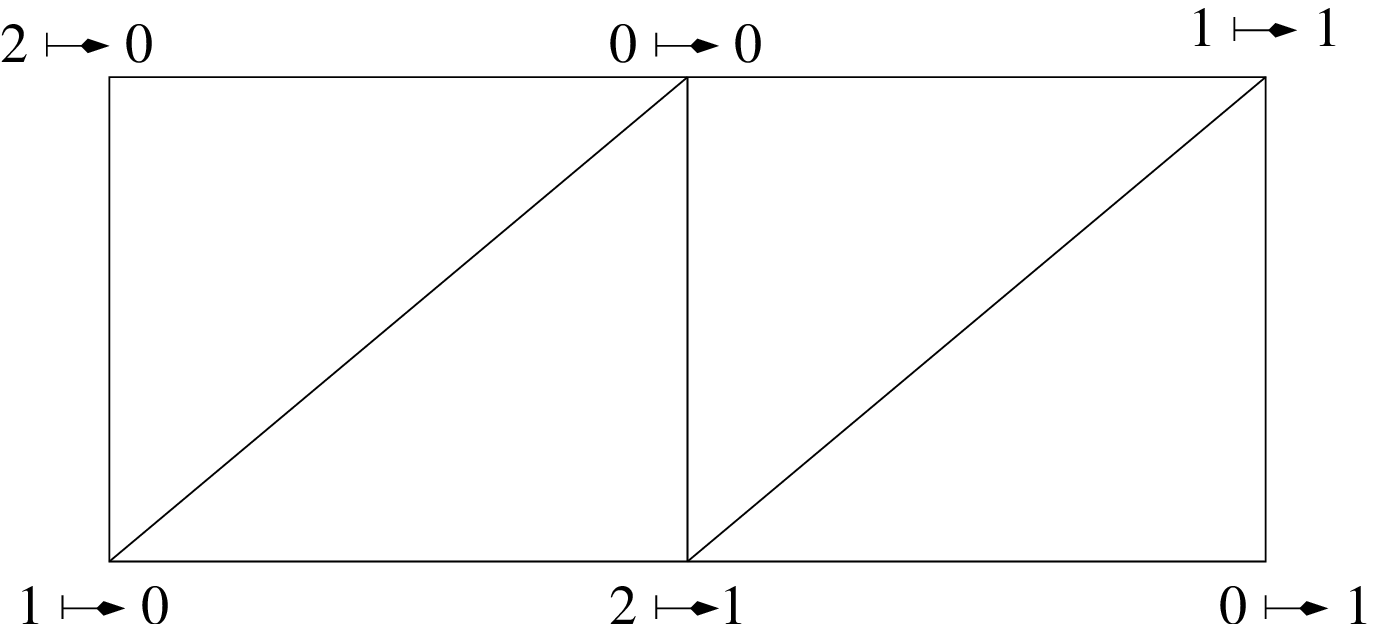}$$
\begin{center}
The simplicial complex $Y_{2 \to 1}$.
\end{center}
\end{figure}

\begin{proof} 
First note that $\G_{n} \downarrow [d]$ is a poset since $\pnno$ is. 
We have that  $(T', f') \leq (T,f)$ if $T' 
\subseteq T$ and $f'$ is the restriction of $f : \G_{n}(T) \to [d]$.
We may identify $T$ with $S = \G_{n}(T) \subset [n]$ and refer only
to the subsets $S$ of $[n]$, to obtain an indexing set for faces as in 
the definition of $Y_{n\to d}$.  That $|\G_{n} \downarrow [d]|$ is 
isomorphic to the subdivision of $Y_{n \to d}$
now follows from the fact that the realization of the face
poset of a complex is isomorphic to its barycentric subdivision.

The vertices of $Y_{n\to d}$
are labeled by morphisms $i \mapsto j$, with $i \in [n]$ and $j \in [d]$.
We may build $Y_{n\to d}$ using these vertices
by stating that some collection spans a (unique) simplex if and only if 
that collection consistently defines a map from a subset of $[n]$ to $[d]$
(see Figure~\ref{F:Ydn}).   
Any such $(S, f)$ admits a map to, and thus labels a face of, some
$([n], \bar{f})$  by letting $\bar{f}(i) = f(s(i))$ where $s(i)$
is the greatest element of  $S$ which is less than or equal to $i$.
Thus, $Y_{n\to d}$ is an $n$-dimensional simplicial complex.   
\end{proof} 

\begin{proof}[Proof of \refT{Gcofin}]
By definition, to show $\G_n$ is left cofinal is to show $|\G_n \downarrow [d]|$,
which by \refL{l1} is homeomorphic to 
$Y_{n\to d}$, is contractible.  

Let $A_{n,d}$ denote the set 
of order-preserving maps from $[n]$ to $[d]$, which by definition
indexes the $n$-simplices of $Y_{n\to d}$.   We put a partial
order on such maps, and thus these $n$-simplices, 
where $f \geq g$ if $f(i) \geq g(i)$ for all $i$.  In this ordering
 $A_{n,d}$ has unique minimal and maximal elements given by $f(i) = 0$ for
all $i$ and $f(i) = d$ for all $i$ respectively.  
Note that any two simplices in $Y_{n \to d}$ which share a face are
comparable in this ordering.

In all cases except for
the minimal element, an $f$ in $A_{n,d}$ shares a face with a smaller $g$,
defined by decreasing a single value of $f$.  If $f(i) = f(i+1)$ for some $i$,
which must be the case when $d < n$, then the face in which $i+1$ is removed
from $S$ (and its value omitted) is not shared by any smaller $g$.  
If there is no $f(i) = f(i+1)$ then $d=n$ and $f$ is the identity, in which case the
face defined by removing $0$ is not shared by any smaller $g$.
Also, all simplices have at least one face not shared by 
smaller simplices.  We may construct a simple homotopy, retracting 
simplices one by one in accordance with our partial ordering.  At each
step we retract the interior of a simplex along with its faces not shared with
lower simplices onto its faces which are shared with lower simplices.

A priori, such a simple homotopy will retract onto the
$(n-1)$ skeleton of $Y_{n \to d}$.
Through closer analysis, we establish
contractibility.
Let $\sigma_{i\mapsto j}$ denote the simplex which minimal among those with
$f(i) = j$, so its other vertices  are labeled  $f(k) = j$ for
$k>i$ and $f(k) = 0$ for $k<i$.
For example, the maximal simplex is $\sigma_{0 \mapsto d}$.
When all simplices greater than $\sigma_{i \mapsto j}$
have been retracted, it will be the only simplex containing the $f(i) = j$ vertex.
So when $\sigma_{i \mapsto j}$ is retracted, we do so onto the face opposite 
from this $f(i) = j$ vertex and then remove it.  By this process we will eventually retract onto the
minimal simplex, establishing contractibility of $Y_{n \to d}$.
\end{proof} 

\medskip 
We end this section with some informal remarks.
The spaces $C'_n[M, \partial]$ are more familiar than $C'_n\la 
M, \partial \ra$ and are manifolds with corners, so it would be preferable
to use them for a cosimplicial model.  One option would be to 
define an $A_\infty$ cosimplicial space, and another option would be to 
enlarge the cosimplicial category.  
We took the latter approach in an early version of this paper. 
But because  the algebra of our spectral sequences is ultimately cosimplicial, we 
opted for using $C'_{n} \la M, \partial \ra$.   We are also in a better
position to take advantage of some recent advances in cosimplicial machinery
\cite{McCl99, Sinh04}.

\medskip

The analogy between our 
cosimplicial model for a knot space and the cosimplicial model for loop spaces 
is made precise in Section~5 of \cite{Sinh04} and Section~5 of \cite{Voli03}.  
Briefly, our machinery can be applied for immersions of an
interval in $M$, namely ${\rm Imm}(I, M)$, as well.  Because immersions may 
self-intersect globally, the $n$th degree approximation 
from embedding calculus is a homotopy limit over $\pnno$ of spaces 
${\rm Imm}(I - \bigcup_{s \in S} I_s, M) \cong (STM)^{\#S -1}$.  Following the 
arguments in this paper we may construct a cosimplicial model which has 
$n$th entry $(STM)^{n}$ and the standard diagonals and projections for
structure maps.  This model is precisely the cosimplicial model for 
$\Omega(STM)$, and this loop space is known to be homotopy equivalent to the 
space of immersions by  theorems of Hirsch and Smale 
\cite{Smal59}.  The spectral sequence in cohomology for this cosimplicial 
model is the Eilenberg-Moore spectral sequence, which is thus analogous 
to the spectral sequences we develop in Section~\ref{S:ss}.

\section{The spectral sequences}\label{S:ss} 

In this section we give spectral sequences which converge to the 
homotopy and cohomology  groups of $\E$, when $M$ is simply connected
and has dimension at least 
four.  We focus particular attention on the case of $M = I^{N+1}$. 

Recall from \cite{BK} that the spectral sequence for the homotopy groups of a cosimplicial 
space $\X^{\bullet}$  is simply the spectral 
sequence for the tower of fibrations  
$ \text{Tot}_0 \X^{\bullet} \leftarrow \text{Tot}_1 \X^{\bullet} \leftarrow \cdots.$ 
Let $s^i$ to denote the codegeneracy
maps of a cosimplicial space and recall that we are using $\delta^{i}$ as the
$i$th face map for our cosimplicial model $\C^\bullet \la M \ra$.

\begin{theorem} \label{T:piss} 
Let $M$ have dimension four or greater.  There is a second quadrant 
spectral sequence converging to $\pi_*(\E)$ whose $E_1$ term is given by
$$E_1^{-p,q} = \bigcap_{i} \; {\text{ker}} \; {s^i}_* \subseteq 
\pi_q\left(C'_p \la M \ra \right). $$ The $d_1$ differential is the restriction to 
this kernel of the map $$\sum (-1)^i {\delta^i}_* \colon 
\pi_q(C'_{p-1} \la M \ra) \to \pi_q(C'_p \la M \ra).$$ 
\end{theorem} 

Before establishing this theorem, we recall some general facts about
homotopy groups of the totalization of a cosimplicial
space.   Looking at the tower of fibrations of the $\text{Tot}_{n}$, 
let $L_n \X^{\bullet}$ be the fiber of the map 
$\text{Tot}_n \X^{\bullet} \leftarrow \text{Tot}_{n-1} \X^{\bullet}.$
If the connectivity of the layers $L_n \X$ tend to infinity, then  $\pi_i(\text{Tot}_n \X^{\bullet})$ 
is independent of  $n$ for $n$ sufficiently large, which implies convergence of the spectral sequence.
As noticed by Goodwillie \cite{Good07}, among others, there is a nice model for 
$L_n \X^{\bullet}$ in terms
of cubical diagrams.  Just as the functor $\G_{n}$ of \refD{gn}
makes a cubical diagram from the cofaces of a cosimplicial diagram, there is a
complementary functor which makes a cubical diagram from the codegeneracy maps.

\begin{definition}
Let ${\G}^{!}_{n} : \pnn \to  \bf{\Delta}_{n}$ be the functor which on objects sends
$S$ to $[n - \#S]$.  On morphisms, it sends the inclusion $S \subset S'$ to the
composite $$[n - \#S] \cong ([n] - S) \overset{p}{\longrightarrow} ([n] - S') \cong [n - \#S'],$$
where the first and last maps are the order-preserving isomorphisms and
where $p$ sends $i \in [n] - S$ to the largest element of $[n] - S'$ which is less than
or equal to $i$.
\end{definition}

Recall from Definition~1.1(b) of 
\cite{GoII} that one way in which the total fiber of a cubical diagram may be defined
is as the homotopy fiber of the canonical map from the initial space in the cube to the homotopy
limit of the rest of the cube.

\begin{theorem}\cite{Good07}\label{T:layer}
The $n$th layer in the Tot tower of a fibrant cosimplicial space $\X^{\bullet}$, 
namely $L_{n} \X^{\bullet}$, is homotopy
equivalent to $\Omega^{n}$ of the total fiber of $i_{n}\X^{\bullet} \circ {\G}^{!}_{n}$.
\end{theorem}

Using this theorem to identify layers in the Tot tower, we can establish the
connectivity results needed for convergence of our homotopy spectral sequence
for knot spaces.

\begin{proof}[Proof of \refT{piss}]
As mentioned above, it suffices to show that the spaces $L_{n} \C^\bullet \la M \ra$ 
have connectivity which grows with $n$.
We apply \refT{layer} and freely use basic facts about cubical 
diagrams from Section~1 of \cite{GoII}.
Let $\D^{!}_{n} \la M \ra$ denote the cubical diagram 
$i_{n} \C^\bullet \la M \ra \circ {\G}^{!}_{n} $.  As mentioned above,
the structure maps of $\D^{!}_{n} \la M \ra$ are given by the  
codegeneracies of $\C^{\bullet} \la M \ra$.

Unraveling definitions, $\D^{!}_{n} \la M \ra$ sends $S$ to  $C'_{n - \#S} \la M, \partial \ra$ 
with structure maps defined by forgetting a point (and associated tangent vector)
from a configuration.  We analyze the total fiber of this cubical diagram through a sequence
of equivalences.  First note that by \refP{proj} and part~\ref{heq} of \refT{mainlist}, we may
replace the spaces  $C'_{i} \la M \ra$ by $C'_{i}(M)$ to obtain an equivalent cubical diagram
$\D^{!}_{n} (M)$.

Next, by naturality of pull-backs there are maps of fiber bundles 
\begin{equation}\label{E:fiber}
\begin{CD}
(S^{m-1})^{n} @>s^i>> (S^{m-1})^{n-1}\\
@VVV @VVV \\
C'_n ( M)    @>s^i>> C'_{n-1} (M ) \\
@VVV @VVV \\
C_n ( M ) @>s^i>> C_{n-1} (M) , 
\end{CD}
\end{equation}
where $m$ is the dimension of $M$ and all horizontal maps $s^i$ are projections.  
These assemble to define
maps of cubical diagrams $$(S^{m-1})^{\pnn!} \to \D^{!}_{n}(M) \to \overline{\D}^{!}_{n}(M),$$
where $\overline{\D}^{!}_{n}(M)$ is obtained from $\D^{!}_{n}(M)$ by applying projections
which forget tangential data at each entry.
Here $(S^{m-1})^{\pnn!}$ sends $S$ to $(S^{m-1})^{\und{n} - S}$, and 
its structure maps are given by projections.    
The total fiber construction for a cubical diagram is functorial, in particular
when using Definition~1.1 of \cite{GoII}.  Moreover, the total
fiber functor sends a fiber 
sequence of cubical diagrams (that is, a sequence of cubical diagrams which
is a fiber sequence at each entry) to a fiber sequence of spaces.  

We claim that the total
fiber of $(S^{m-1})^{\pnn!}$ is contractible for $n\geq 2$, which would then imply that
the total fibers of $\D^{!}_{n}(M)$ and $\overline{\D}^{!}_{n}(M)$ are equivalent.
The easiest way to see this contractibility 
is to use an alternate inductive definition of total fiber.  
The total fiber an $n$-cube $Z$ is homotopy equivalent to that of the $(n-1)$-cube
$\phi Z$, which sends ${S}$ to the homotopy fiber of the map from $Z_{S \cup n}$
to $Z_{S}$.  In our case where $Z = (S^{m-1})^{\pnn!}$, the homotopy fiber of  
$Z_{S \cup n} \to Z_{S}$ is always simply
$S^{m-1}$ and the structure maps in $\phi Z$ are all identity maps.
Applying $\phi$ again, $\phi^{2}(S^{m-1})^{\pnn!}$ 
has a point in every entry, yielding a contractible total fiber.  

We may now focus on $\overline{\D}^{!}_{n}(M)$, and
consider $\phi \overline{\D}^{!}_{n}(M)$.  Since the structure maps of $\overline{\D}^{!}_{n}(M)$,
 namely projections of configuration spaces, are fibrations we take
 fibers to see that 
 $\phi \overline{\D}^{!}_{n}(M)$ sends $S$ to  $M$ 
 with $(n-1 - \#S)$ points removed.  
 In order for structure maps to commute,  
 we need to fix an embedding $f$ of $\und{n-1}$ in $M$, so that the $S$ entry of 
 $\phi \overline{\D}^{!}_{n}(M)$ is $M - f(\und{n-1} - S)$.  The structure maps are then inclusions.

 We analyze the  total fiber of $\phi \overline{\D}^{!}_{n}(M)$  by
 applying the Blakers-Massey theorem.  Inclusions of codimension zero
 open submanifolds are cofibrations, and $\phi \overline{\D}^{!}_{n}(M)$
 is a push-out cube, so it is strongly-coCartesian.  Moreover, each initial inclusion
 $M -  f(\und{n-1}) \hookrightarrow M - f(\und{n-1} -  i)$ is an $(m-1)$-connected map.
  Indeed, this inclusion is
 surjective on $\pi_{d}$ for $d \leq m-1$ since by basic transversality any map from $S^{d}$
 to $M - f(\und{n} - i)$ is homotopic to one which avoids the point $i$ and thus lies in
 $M - f(\und{n})$.  Applying the same argument to homotopies gives that this inclusion
 is injective on $\pi_{d}$ for $d < m-1$.
 The Blakers-Massey theorem (stated as Theorem~2.3 of \cite{GoII})
 applies to give that the total fiber of $\phi \overline{\D}^{!}_{n}(M)$
 is $(n-1)(m-2) + 1$-connected.
 
 Because we have 
 equated the total fiber of $\phi \overline{\D}^{!}_{n}(M)$ with that of $\D^{!}_{n} \la M \ra$,
 this shows that the connectivity of $L_{n} C^{\bullet} \la M \ra$ is $(n-1)(m-2) + 1 - n = (n-1)(m-3)$ 
 (where the $-n$ appears because of the  $\Omega^{n}$ in \refT{layer}), 
 which does increase with $n$ as needed
 when $M$ has dimension four or greater.
 \end{proof}
 
 Our connectivity estimates in proof give the following.
 
 \begin{corollary}
The group $E_{1}^{-p, q}$ in the spectral sequence of \refT{piss} vanishes if
$q \leq (p-1)(m-2)$, where $m$ is the dimension of $M$.
\end{corollary}

We give more explicit computations of the rational
homotopy spectral sequence when 
$M$ is $I^{N+1}$ or equivalently $\R^{N} \times I$ in \cite{2}.  
The rows of this spectral sequence
were conjectured by Kontsevich in \cite{Kont00}.   
Lambrechts and Tourtchine reformulate this spectral
sequence in terms of graph homology in \cite{LaTo06}.

\medskip 

We now discuss the cohomology spectral sequence.   We proceed by 
taking the homology spectral sequence first studied in  
\cite{Re70} and dualizing through the universal coefficient theorem. 
The construction and convergence of the homology spectral sequence is more delicate  
than that of the homotopy spectral sequence \cite{Bott77, Bous86, Ship96},
but we will see that convergence follows from similar connectivity estimates.  

\begin{theorem} \label{T:cohss} 
Let $M$ be a simply connected manifold with two distinguished
points on its boundary of dimension four or greater.  There is 
a second quadrant spectral sequence converging to $H^*(\E; \Z/p)$ whose 
$E_1$ term is given by
$$E_1^{-p,q} = {\rm{coker}} \; \sum (s^i)^* \colon H^q(C'_{p-1}(M); 
\Z/p) \to H^q(C'_p(M); \Z/p).$$ The $d_1$ differential is the passage to 
this cokernel of the map $$\sum_{i=1}^{p+1} (-1)^i (\delta^i)^* \colon H^q(C'_p(M); 
\Z/p) 
\to H^q(C'_{p-1}(M); \Z/p).$$ 
\end{theorem} 

This is sometimes called
the normalized $E_1$ term, quasi-isomorphic to one with $E_1^{-p, q} = H^q(C'_p(M); \Z/p)$
and $d_1$ as given but without passage to cokernels.

The case of knots in Euclidean spaces serves as both a central example and a base case
for the proof of this theorem.  We now analyze these $E_{1}$ groups in the 
setting of  $M = I^{N+1}$, where there is  a more combinatorial
description.  First note that $C_n \la I^{N+1}, \partial \ra$ is homotopy equivalent to 
$C_n(\R^{N+1})$.  Recall from \cite{FCoh95, FCoh76, FCoh75.2, Arno68} that 
$H^*(C_n(\R^{N+1}))$
is generated by classes $a_{ij}$, which are pulled back by $\pi_{ij}^*$
from the generator of $H^N(S^N)$, with $1 \leq i \neq j \leq n$.
In addition to graded commutativity there are relations 
\begin{gather} 
{a_{ij}}^2 = 0  \label{basics1} \\ a_{ij} = (-1)^{N+1} a_{ji}\\ 
a_{ij} a_{jk} + a_{jk} a_{ki} + a_{ki} a_{ij} = 0 \label{Jacobi}.
\end{gather} 
By the K\"unneth  theorem, the cohomology ring of $C_n' (\R^{N+1})$ 
is isomorphic to that of $C_n(\R^{N+1})$ tensored 
with an exterior algebra on $n$ generators, which we call $b_1$ 
through $b_n$.

There is a natural
description of these cohomology groups in terms of 
graphs  with vertex set $\n$ and with edges which are 
oriented and ordered.  Let  $\Gamma_{n}$
denote the free module generated by such graphs, which is a ring by taking
the union of edges of two graphs in order to multiply them.  Then $\Gamma_n$
surjects onto $H^*(C_n'(\R^{N+1}))$ by sending a generator, which is a
graph with a single edge from $i$ to $j$, to $a_{ij}$ if $i\neq j$ or $b_i$ if $i=j$.
Let $\iop'(n)$ (respectively $\iop'_{o}(n)$ if $N$ is odd)
denote $\Gamma_n$ modulo the kernel of this homomorphism,
which by construction is isomorphic to $H^*(C_n'(\R^{N+1}))$.   

Because the cohomology 
of $C_n'(\R^{N+1})$ is torsion free we work integrally in  the cohomology spectral sequence
for $C^\bullet \la I^{N+1} \ra$.   The map 
$s^k : C_n'\la I^{N+1}, \partial \ra \to C_{n-1}' \la I^{N+1}, \partial \ra$
restricts to the map on open configuration spaces which forgets the
$k$th point in a configuration, whose effect on cohomology we now determine.

Let $\sigma_{\ell}: \und{n-1} \to \und{n}$ be the order-preserving inclusion for which
$\ell$ is not in the image. For $i,j \neq \ell$, the map $\pi_{ij} \colon C_n(I^{N+1}) \to 
S^N$ factors as 
$\pi_{\sigma_\ell(i) \sigma_\ell(j)} \circ s^\ell$.  Hence
$a_{ij} \in H^N(C_{n-1}(I^{N+1}))$ maps to $a_{\sigma_k(i) \sigma_k(j)} 
\in H^N(C_n(I^{N+1}))$ under $s_k^*$.  Translating through the 
isomorphism with the modules $\iop(n)$, 
$(s^\ell)^*$ takes a graph with $n-1$ vertices, relabels the 
vertices according to $\sigma_{\ell}$ and adds a vertex, not attached to any 
edges nor marked, labeled $\ell$.  Hence the sub-module generated by the 
images of $(s^\ell)^*$ is the sub-module $D_n$ of graphs in
which at least one vertex is not attached to an edge.  The quotient map 
$\iop'_m(n) \to \iop'_m(n) / D_n$ is split, so that $\iop'_m(n) / D_n$ is 
isomorphic to the submodule  of $\iop'_m(n)$ generated by graphs in 
which every vertex is attached to an edge, which we call $\overline{\iop'_m(n)}$. 

Next we identify the homomorphisms $(\delta^\ell)^*$,
where $\delta^\ell$ ``doubles'' the $\ell$th point in a configuration.  
Let $\tau_\ell$ (by slight abuse, after \refD{cofaces}) be the order-preserving
surjection from $\und{n+1}$ to $\und{n}$ for which $\ell$ and $\ell+ 1$ map to $\ell$.
If $(i,j) \neq (\ell, \ell+1)$ then the composite  
$\pi_{ij} \circ \delta^\ell$ coincides with $\pi_{\tau_\ell(i), 
\tau_\ell(j)}$, which implies that $a_{ij}$ maps to  
$a_{\tau_\ell(i), \tau_\ell(j)}$ under $(\delta^\ell)^*$.  On the
other hand, $\pi_{\ell, \ell+1} \circ \delta^\ell$ is  
the projection of $C'_n(I^{N+1})$ onto the $\ell$th tangential factor of $S^N$,  
so that $(\delta^\ell)^*(a_{\ell, \ell+1}) = b_\ell$.  We 
extend these computations to define $(\delta^\ell)^*$ on all 
of $H^*(C'_n(I^{N+1}))$ using the cup product.   
In terms of graphs, let $c_\ell$ be the map on $\iop'(n)$ defined 
by identifying the $\ell$th and 
$\ell+1$st vertices in a graph and relabeling.    See Figure~\ref{F:chord} for an illustration,
which will be familiar to those who have read \cite{Vass}.

\skiphome{
\begin{figure}\label{F:chord}
$$\includegraphics[width=14cm]{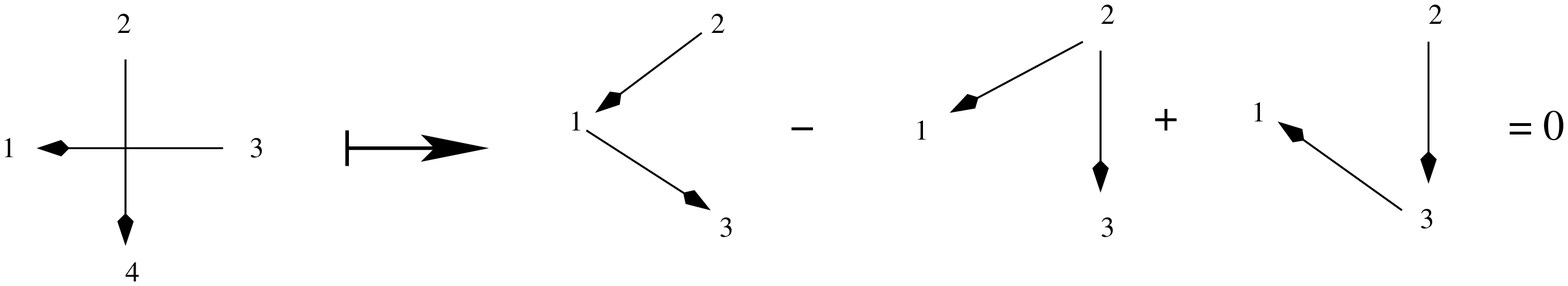}$$
\begin{center}
Figure~\ref{F:chord} - An illustration of $d_1 = \sum (-1)^{i} {\delta^{i}}^{*}$.
\end{center}
\end{figure}
}

Let $\overline{\iop'_m(n)} \subset \overline{\iop'(n)}$ denote the submodule generated by
graphs with $m$ edges.   Our graphical analysis of the maps $s^{i}$ and $d^{i}$
above leads to the following.

\begin{proposition}
The group named in \refT{cohss} as $E_{1}^{-p, q}$ for the cohomology spectral
sequence for the space of knots in $M = I^{N+1}$, namely 
 ${\rm{coker}} \; \sum (s^i)^* \colon H^q(C'_{p-1}(I^{N+1}); 
\Z/p) \to H^q(C'_p(I^{N+1}); \Z/p)$, is isomorphic to 
$\overline{\iop'_m(p)}$ for $q = mN$ and is zero otherwise.  Moreover, 
$d_1 = \sum (-1)^\ell c_\ell$
\end{proposition}

We obtain a vanishing result  needed for convergence
of the cohomology spectral sequence. 

\begin{corollary}\label{C:vanish} 
The group named in \refT{cohss} as $E_1^{-p,q}$ for the cohomology spectral sequence
for knots in $I^{N+1}$
vanishes when $q < \frac{N}{2}  p$ and when $q > N (2p -1)$ . 
\end{corollary}

\begin{proof}

For the lower vanishing line we argue combinatorially.  The 
module $\overline{\iop'_m(p)}$ is zero if $m< p/2$, as each edge connects at most two vertices.

For the upper vanishing line we argue topologically.  The
Leray-Serre spectral  sequences for the fiber bundles
$(I^{N+1} - i \; {\text{points}}) \to C_{i+1}(I^{N+1}) \to C_i(I^{N+1})$
yield that the cohomology of $C_p(I^{N+1}) $ 
vanishes above degree $N(p-1)$.  Thus the cohomology of $C_p(\R^{N+1})
\times (S^k)^p$ vanishes above degree $N(2p-1)$.  
\end{proof}

Tourtchine improves the upper vanishing line and makes computations of the groups
adjacent to that vanishing line in \cite{Tour04}.

The lower vanishing line is the basis for establishing convergence of this spectral sequence,
which we do now.

\begin{proof}[Proof of \refT{cohss}]
For a cosimplicial space $\X^\bullet$ let $\overline{H}_q(\X^p)$, which
we call the normalized homology of $\X^p$, be the intersection of the kernels of the 
codegeneracy maps $s^i \colon H_q(\X^p) \to H_q(\X^{p-1})$.  Theorem 3.4 of \cite{Bous86} states 
that the mod-$p$ homology spectral sequence of a cosimplicial space
$\X^\bullet$ converges (strongly) when three conditions are met, namely 
that $\X^p$ is simply connected for all $p$, 
$\overline{H}_q(\X^p) = 0$ for $q\leq p$, and for any given $k$ only 
finitely many $\overline{H}_q(\X^p)$ with $q-p = k$ are non-zero. 
The last two conditions are satisfied if $\overline{H}_q (\X^p)$
vanishes for $q < cp$  for some $c>1$, which we call the 
{\em{vanishing condition}}.  

The vanishing condition for 
$\C^\bullet \la I^{N+1} \ra$ is
given as the lower vanishing line in Corollary~\ref{C:vanish} (after applying
the universal coefficient theorem) with $c = \frac{N}{2}$, 
yielding convergence in the Euclidean case in dimensions four and higher.

For general simply connected $M$ of dimension $m \geq 4$,
we analyze the degeneracy maps 
$C'_p \la M, \partial  \ra \to C'_{p-1} \la M, \partial \ra$ through the
fiber sequences $(S^{m-1})^{\und{p}} \to C'_{p}(M) \to C_{p}(M)$.
As noted in the proof of  \refT{piss}, these fiber sequences commute with
degeneracy maps.
Because the Leray-Serre spectral sequence is natural, it 
suffices to check the vanishing condition for the normalized homology of the
fiber and base separately.  The vanishing condition for the fiber
is a standard feature of the Eilenberg-Moore spectral sequence
for the homology of $\Omega S^{m-1}$, where the degeneracy maps
are projections $(S^{m-1})^{i} \to (S^{m-1})^{i-1}$ and the vanishing
constant is $c = m-1$.

To check the vanishing condition for the normalized homology of  $C_p(M)$,
we follow Totaro  \cite{Tota96} and study 
the Leray spectral sequence of the inclusion $i : C_p(M) \hookrightarrow M^{\underline{p}}$
(an equivalent spectral sequence was first constructed by Cohen and Taylor \cite{FCoh78}). 
The squares
$$
\begin{CD}
C_{p}(M) @>s^{j}>>   C_{p-1}(M) \\
@ViVV                    @ViVV\\
M^{\und{p}}  @>s^{j}>>  M^{\und{p-1}}
\end{CD}
$$
commute, where the second horizontal arrow
is the obvious projection map, naturally labeled $s^{j}$ because it serves
as a codgeneracy map in the cosimplicial model for $\Omega M$.
We use naturality of the Leray spectral sequence to  check 
the vanishing conditions on the cohomology of the base and the stalks
separately.  

On the base $M^{\und{p}}$, the vanishing condition is identified with
that for  the Eilenberg-Moore spectral sequence for $\Omega M$, which  is well known
for simply connected $M$.  
By the K\"unneth theorem the normalized homology of $M^{\und{p}}$ vanishes
below degree $(k+1)p$ where $k$ is the connectivity of $M$.  

The stalk of $i$
over  a point in $M^{\und{p}}$ is the cohomology of 
the product $\prod_{S} C_{\#S}(\R^{m})$, where $S$ ranges of the components
of the partition of $\und{p}$ defined by $i \sim j$ if $x_{i} = x_{j} \in M$.
The homology of this product is naturally a submodule of $\iop'(p)$, namely
that spanned by graphs whose edges may only connect elements of the
same component of the partition.  Moreover, the image of the maps 
${s^{j}}^{*}$ are still graphs with some vertex not connected to any edge.
\refC{vanish} applies to establish the vanishing condition for these
stalks with the same constant as for knots in Euclidean spaces, namely
$\frac{m-1}{2}$.  

Using these spectral sequences to combine 
the vanishing estimates for the normalized homology 
of $(S^{m-1})^{\und{p}}$, $M^{\und{p}}$, and the stalks of 
$i : C_p(M) \hookrightarrow M^{\und{p}}$, we establish
the vanishing condition for the normalized homology of
$C'_p (M)$.  The constant is the lesser of $\frac{m-1}{2}$
and $k+1$, where $m$ is the dimension of $M$ and $k$ is its connectivity.
\end{proof}

In the course of proof, we have established a lower vanishing line for our
cohomology spectral sequence for $\E$ which coincides
with either that for  $\Omega M$ or that for ${\rm Emb}(I, I^{m})$.

\begin{corollary}
The group ${E}_{1}^{-p, q}$ in the cohomology spectral sequence
for ${\rm Emb}(I, M)$ of \refT{cohss} vanishes when $q$ is smaller
than the lesser of $\frac{m-1}{2} p$ and $(k+1) p$, 
where $m$ is the dimension of $M$ and $k$ is its connectivity.
\end{corollary}

\end{document}